\documentclass[11pt]{amsart}
\usepackage{amsmath,amssymb}
\numberwithin{equation}{section}

\newtheorem{itheorem}{Theorem}
\newtheorem{theorem}[equation]{Theorem}

\newtheorem{lemma}[equation]{Lemma}
\newtheorem{corollary}[equation]{Corollary}

\theoremstyle{definition}
\newtheorem{definition}[equation]{Definition}
\newtheorem{remark}[equation]{Remark}
\newtheorem{example}[equation]{Example}

\def\C{\mathbb C}
\def\N{\mathbb N}

\def\R{\mathbb R}
\def\Z{\mathbb Z}

\def\S{\mathcal{S}}
\def\H{\mathcal{H}}
\def\U{\mathbf{U}}
\def\V{\mathcal{V}}
\def\rational{\mathcal{R}}

\def\st{\,|\,}
\DeclareMathOperator{\LinSpan}{span}

\DeclareMathOperator{\spec}{spec}
\begin{document}

\title{Hecke operators on rational functions}
\author{Juan B. Gil}
\address{3000 Ivyside Park, Penn State Altoona\\ Altoona, PA 16601}
\email{jgil@psu.edu}
\author{Sinai Robins}
\thanks{The second author was supported by the NSA Young Investigator Award
{\tiny MSPR-00Y-196}}
\address{Department of Mathematics\\ Temple University\\
 Philadelphia, PA 19122}
\email{srobins@math.temple.edu}

\keywords{Rational functions, Hecke operators}
%\subjclass[2000]{Primary ; Secondary }
%\begin{abstract}
%\end{abstract}

\maketitle
\tableofcontents
%%%%%%%%%%%%%%%%%%%%%%%%%%%%%%%%%%%%%%%%%%%%%%%%%%%%%%%%%%%%%%%%%%%%
\section{Introduction}
%%%%%%%%%%%%%%%%%%%%%%%%%%%%%%%%%%%%%%%%%%%%%%%%%%%%%%%%%%%%%%%%%%%%

We study the vector space of all linear recurrence sequences over
the reals, by defining linear operators that sift out arithmetic
progressions from the sequence.  We call these linear operators
Hecke operators, by analogy with the theory of automorphic forms,
and we develop their spectral theory completely. Because the
generating function of any linear recurrence sequence is a
rational function (with nonzero constant term in the denominator)
this study is equivalent to the action of our linear operators on
rational functions.  Although we borrow terminology from the
traditional theory of Hecke operators on modular forms, prior
knowledge of Hecke operators is not assumed here, since both the
problems and the methods herein are grounded in the new context of
rational functions.

To begin our study of Hecke operators, we let $\rational$ be the
vector space of all rational functions $f(x)=A(x)/B(x)$ with real
coefficients such that $\deg A(x) < \deg B(x)$, and such that
$B(0)\not= 0$.  We note that $A$ and $B$ are not restricted to be
relatively prime, and that the degree restriction is only assumed
for ease of notation. Given a rational function $f\in \rational$
whose Taylor series is $f(x) = \sum_{n=0}^\infty a_n x^n$ and
given a positive integer $p\in\N$, we define the \emph{Hecke
operator} $U_p:\rational\to \rational$ by
\begin{equation}\label{HeckeOpDefinition}
U_p f(x) = \sum_{n=0}^\infty a_{pn} x^n.
\end{equation}

It turns out that the class of rational functions which are
eigenfunctions of at least one of the Hecke operators defined above
generate the subspace $\rational_{qp} \subset \rational$ of all
rational functions with poles at the roots of unity.

An equivalent description of this class of rational functions can be
given by noting that each rational $f = \sum_{n=0}^\infty a_n x^n \in
\rational_{qp}$ has coefficients that are \emph{quasipolynomials} in
$n$ (by the standard Theorem~\ref{RationalFunctions} below), and hence
our use of the subscripts {\it qp} in $\rational_{qp}$.  That is,
\[ a_n = c_d(n) n^d + c_{d-1}(n) n^{d-1} + \dots + c_0(n), \]
where each $c_j(n) \in \mathbb{Q}$ is a periodic function on $\mathbb{Z}$.
There is a large body of knowledge on Ehrhart quasi-polynomials whose
generating functions give rational functions in $\rational_{qp}$.
Thus $\rational_{qp}$ provides a rich source of functions that arise
naturally in the theory of lattice point enumeration in rational
polytopes and combinatorial geometry (see for example
\cite{MR2000k:52014}, \cite{MR98e:11117a}, and \cite{Stan97}).  We
note that the reader does not, however, require any previous knowledge
in this field for the analysis presented here.

In Section~\ref{sec_Spectrum} we study the spectral properties of the
Hecke operators on rational functions and show that they have discrete
spectra. The first result concerning the structure of eigenfunctions
is the following:

\begin{itheorem}[Involution Property]
If $f(x)=A(x)/B(x)$ is an eigenfunction of $U_p$, then all roots
of $B(x)$ are roots of unity and we have the identity
\begin{equation*}
 x^d B(\tfrac{1}{x}) = (-1)^d B(x).
\end{equation*}
Moreover, if $U_p f = \lambda f$ with $\lambda \not= 1$, then
\[ U_p (f(\tfrac{1}{x})) = \lambda f(\tfrac{1}{x}).  \]
That is,  $f(\tfrac{1}{x})$ is another eigenfunction of $U_p$ with
the same eigenvalue $\lambda$, and with the same denominator
$B(x)$.
\end{itheorem}

The involution $x \mapsto \frac{1}{x}$ on eigenfunctions plays an analogous
role to the Fricke-Atkin-Lehner involutions on eigenforms.  Here we uncover
more properties of eigenfunctions, reducing the problem of computing
eigenvalues of an infinite dimensional linear operator to
the problem of computing eigenvalues of a finite matrix.

The main result of this section is described by:

\begin{itheorem}[The Spectrum]
Let $p$ be any positive integer greater than $1$. Then the point
spectrum of $U_p$ on the vector space $\rational$ is
\[ \spec(U_p) = \{\pm p^{k} \st k\in\N\}\cup \{0\}. \]
\end{itheorem}

The eigenfunctions of any Hecke operator obey a rigid structure
theorem, described in Section~\ref{sec_Structure}, that makes them
appealing and easy to work with.  There is far more structure in
$\rational$ due to these eigenfunctions than has been hitherto apparent.

The main structure theorem is the following:
\begin{itheorem}[Structure Theorem]
Let $f(x)=A(x)/B(x)=\sum a_n x^n$ be an eigenfunction of $U_p$ for
some integer $p>1$, associated to an eigenvalue $\lambda_p
\not=0$. If $B(x)= \prod_{j=1}^d (1-\gamma_j x)$, then there is an
integer $\kappa$ dividing the degree $d$, and an integer $L$ such that
\begin{equation*}
 a_n = n^{\kappa-1}\sum_{j=1}^{d/\kappa} C_{j}\,
 e^{\frac{2\pi i  \ell_j}{L}n}, \text{ for all} \  n \geq 0,
\end{equation*}
where each pole of $f$ is given by $\gamma_j =e^{\frac{2\pi i
\ell_j}{L}}$, $\ell_j\in\N$, and the constants $C_{j}\in \C$ are determined
by the initial conditions of the linear recurrence sequence
$\{ a_n \}$. We note that each pole $\gamma_j$ must occur with the
same multiplicity $\kappa$.
\end{itheorem}

When we consider rational functions that are simultaneous
eigenfunctions of a family of Hecke operators, we discover that there
is a natural character $\chi_f \pmod{L}$ that comes into the spectrum.
We also get the curious phenomenon of ``partial characters'' whenever we
find a rational function $f$ that is an eigenfunction of some, but
\emph{not all} Hecke operators.  

Here $L$ is called the \emph{level} of $f$, and is defined for any $f
\in \rational_{qp}$ as the least common multiple of the orders of all
the roots of unity that comprise the poles of $f$. For example, if the
poles of $f$ are $e^{2\pi i/4}$ and $e^{2\pi i / 5}$, then $f$ has
level $L=20$.  There is also a notion of the \emph{weight} $\kappa$ of
an eigenfunction $f$, arising naturally in the structure theorem
above, and defined by the common multiplicity of the poles.

In Section~\ref{sec_PGR} we decompose the infinite dimensional
vector space $\rational_{pq}$ into finite dimensional subspaces,
by using the weight and level of an eigenfunction as the grading
parameters.  We let $\V_{\kappa,L}(U_3, U_5)$, for example, denote
the finite dimensional vector space of eigenfunctions of (at
least) the Hecke operators $U_3$ and $U_5$ that have weight
$\kappa$ and level $L$.

By a further analogy with modular forms, the set of simultaneous
eigenfunctions  for the full Hecke algebra $\mathfrak{H}$ (that
are not in the kernel of any $U_p$) is of special interest. In
Section~\ref{sec_Simultaneous} we give a complete description of
the vector space $\V$ spanned by all of the simultaneous
eigenfunctions of $\mathfrak{H}$.  The next two results handle the
two separate cases when the level is $L=1$, and $L>1$ for
simultaneous eigenfunctions.

\begin{itheorem}[Simultaneous Eigenfunctions]
Let $f$ be a simultaneous eigenfunction of $\mathfrak{H}$ such
that $f$ is not in the kernel of $U_p$ for any $p$, and
level$(f)=L$. That is, let $f$ be a rational function with the
property that for \textit{every} $p$ there is a $\lambda_p\not=0$
such that $U_p f=\lambda_p f$. Then $L =1$, and
\[ f(x) = C (x\partial_x)^k \left(\frac1{1-x}\right)\]
for some $k\in\N$ and $C\in\C$.
\end{itheorem}

In Corollary~\ref{Character}, the following counterpart to this
result is given for all simultaneous eigenfunctions that have some
of its Taylor coefficients equal to zero.

\begin{itheorem}
Let $L>1$ be a given integer.  Suppose $f(x)=\sum_{n=0}^\infty a_n
x^n$ is a real rational function of level $L$ with $a_0=0$, and
$f(x)$ is a simultaneous eigenfunction of the operators $U_2, U_3,
\dots, U_L$ (i.e. $\;U_m f = \chi_f(m) m^{\kappa -1} f$ for every
$m=2,\dots,L$).

Then $\chi_f$ is the real quadratic character mod $L$, $f$ is in
fact a simultaneous eigenfunction of all the Hecke operators
$U_m$, and in addition we must have
\[  f(x) = a_1 \sum_{n=0}^\infty  \chi_f(n) n^{\kappa -1} x^n. \]
It is worthwhile noting that $f$ can also be written as
\begin{equation*}
  f(x) = a_1 \ (x \partial_x)^{\kappa-1}
  \bigg(\frac{\sum_{j=1}^{L-1} \chi_f(j) x^j}{1-x^L} \bigg).
\end{equation*}
Under the same hypothesis, except with $a_0\not=0$, we conclude
that $\chi_f$ is the principal character and
$f(x)=\frac{a_0}{1-x}$.
\end{itheorem}

Here the differential operator $x\partial_x$ plays the role of the
``weight-raising'' operator in modular forms, because it takes
eigenfunctions of weight $\kappa$ to eigenfunctions of weight $\kappa+1$.

As a curious application of the explicit characterization of $\V$,
the vector space of simultaneous eigenfunctions, we can realize
any finite Euler product of the Riemann zeta function in
Section~\ref{sec_TensorProduct} as the spectral zeta function of a
very explicit operator $\U_S$. More precisely, for any finite set
of primes $S = \{ p_1, \dots, p_n \}$, we define a corresponding
operator $\U_S$ as a finite tensor product $U_{p_1}\otimes\cdots
\otimes U_{p_n}$.  This operator $\U_S$ acts on tensor products of
eigenfunctions, and it turns out that we retrieve any finite piece
of the Euler product for the Riemann zeta function, precisely as
the spectral zeta function $ \zeta_{\U_S}(s) = \sum_{\lambda \in
\spec(\U_S)} \frac{1}{\lambda^s}$ of the operator $\U_S$.

\begin{itheorem}[Euler product]
\begin{equation*}
 \zeta_{\U_S}(s)
 =\zeta_{U_{p_1}}(s)\cdots \zeta_{U_{p_n}}(s)
 =\prod_{p\in S} \frac1{1-p^{-s}}.
\end{equation*}
\end{itheorem}

To extend these ideas to infinite Euler products, we now define
$\H^\infty$ to be the space of products $\mathbf{f}= f_1\otimes
f_2\otimes \cdots$, where $\{f_n\}_{n\in\N}$ is an infinite
sequence of rational functions with the following properties:
\begin{enumerate}
\item There is a finite set $I\subset \N$ such that $f_j\in \V$
for every $j\in I$. \item $f_j= 1$ for every $j\in \N\setminus I$.
\end{enumerate}

For $\mathbf{f}\in \H^\infty$ we define the operator $\U$ by
\begin{equation}\label{BigTensor}
  \U\mathbf{f} =
  (U_{p_{  i_1}} f_{i_1}) \otimes\cdots\otimes (U_{p_{i_{m}}}f_{i_m}),
\end{equation}
where $I=\{i_1,\dots,i_m\}$ is the finite set of positive integers
associated to $\mathbf{f}$, and where $p_{i_k}$ is the $i_k$'th
prime number. Notice that similarly to $\U_S$, the operator $\U$
maps tensor products of rational functions into rational functions
in several variables.

\begin{itheorem}[Riemann zeta function]
The spectral zeta function of the operator $\U$ on $\H^\infty$
satisfies
\begin{equation*}
 \zeta_{\U}(s) = \zeta(s),
\end{equation*}
where $\zeta(s)$ is the Riemann zeta function.
\end{itheorem}
Finally, we conclude with an Appendix that displays explicit
examples of eigenfunctions, to give the reader a better feeling
for the eigenfunctions and eigenspaces that arise.  It is a highly
non-trivial problem to compute the dimensions of the various
vector spaces of eigenfunctions defined by fixing the weight and
the level of admissible eigenfunctions.  Indeed, even computing
dim$(\V_{1,L}(U_p))$, for example, involves the Artin conjecture
for $p$ a primitive root (mod $L$) for infinitely many
integers $L$.

The examples in the Appendix suggest the following interesting
unimodality conjecture about the numerators of all eigenfunctions for
any $U_p$.  More precisely, we can always bring an eigenfunction $f(x)
= A(x)/B(x)$ into a unique canonical form with a denominator $\prod
(1-x^{m_j})$   by multiplying $A(x)$ and $B(x)$ by a common factor if
necessary.  We now assume that $f$ has this form.

\medskip
\noindent
{\bf Conjecture.}
\noindent
Let $f(x) = A(x)/B(x)$ be an eigenfunction of at least one Hecke operator,
written in the unique canonical form given above.
Then the absolute value of the nonzero coefficients of $A(x)$ form a unimodal
sequence.

%%%%%%%%%%%%%%%%%%%%%%%%%%%%%%%%%%%%%%%%%%%%%%%%%%%%%%%%%%%%%%%%%%%%
\section{Some preliminaries}\label{sec_Definitions}
%%%%%%%%%%%%%%%%%%%%%%%%%%%%%%%%%%%%%%%%%%%%%%%%%%%%%%%%%%%%%%%%%%%%

We recall some standard facts about linear recurrence
sequences and their generating functions. The following theorem gives a
characterization of linear recurrence sequences in terms of rational
functions, and gives a closed form for their Taylor
coefficients.  For a proof see  R. Stanley's book \cite[Chapter 4]{Stan97}.

\begin{theorem}\label{RationalFunctions}
Let $\alpha_1, \alpha_2,\dots,\alpha_d$ be a sequence of complex numbers,
$d\in \N$ and $\alpha_d\not=0$. Consider the formal power series
$\sum_{n=0}^\infty a_n x^n$.  The following conditions on the
coefficients $a_n$ are equivalent:
\begin{enumerate}
\item[(i)]
\[ \sum_{n=0}^\infty a_n x^n = \frac{A(x)}{B(x)}, \]
where $B(x)=1+\alpha_1 x + \cdots + \alpha_d x^d$ and $A(x)$ is a
polynomial in $x$ of degree less than $d$.
\item[(ii)]
For all $n\in\N$,
\[ a_{n+d} = -\alpha_1 a_{n+d-1} - \cdots - \alpha_d a_{n}. \]
\item[(iii)]
For all $n\in\N$,
\begin{equation*}
  a_n = \sum_{j=1}^d C_j n^{m_j-1} \gamma_j^n,
\end{equation*}
where each $C_j\in \C$, $m_j$ is a positive integer, and
\[ 1+\alpha_1 x + \cdots + \alpha_d x^d =
   \prod_{j=1}^d (1-\gamma_j x). \]
Each $m_j$ is the multiplicity of the root $\gamma_j$.
\end{enumerate}
\end{theorem}

When all the poles $\gamma_j$ above are roots of unity,
$a_n$ is known as a \emph{quasi-polynomial} in $n$, cf. \cite{Stan97}.

One of our primary goals is to study the spectral properties of
Hecke operators acting on the vector space of rational functions
$\rational$. To this end we must first justify the definition of
the Hecke operator \eqref{HeckeOpDefinition} with a lemma.  That
is, we do not yet know that the image of $U_p f$ is indeed a
rational function in $\rational$.

%%%%%%%%%%%%%%%%%%%%%%%%%%%%%%%%%%%%%%%%%%%%%%%%%%%%%%%%%%%%%%%%%%%%
\begin{lemma}
Given a rational function $f\in \rational$, $U_p f$ is again in
$\rational$.  Moreover, there is a simple algorithm that
constructs the rational function $U_p f$ from the roots of $f$. If
the pole set of $f$ is $\{ \gamma_1, \dots, \gamma_d \}$, then the
pole set of $U_p f$ is $\{\gamma_1^p,  \dots, \gamma_d^p \}$.
\end{lemma}
\begin{proof}
We employ the structure Theorem~\ref{RationalFunctions} to write
the Taylor coefficients of $U_p f$ as
\begin{equation*}
  a_{pn} = \sum_{j=1}^d C_j (pn)^{m_j-1} (\gamma_j^p)^n.
\end{equation*}
Thus, the defining characteristic polynomial for the sought-after
linear recurrence given by $U_p f$ is
\begin{equation}\label{PolyForTp}
\prod_{j=1}^d (1-\gamma_j^p x) =  1+\beta_1 x + \cdots + \beta_d x^d.
\end{equation}
It is clear that the coefficients of this polynomial are real, since the
Taylor coefficients of $f$ are real.  We provide a simple algorithm for
finding the polynomial \eqref{PolyForTp}.

By the fundamental theorem on symmetric polynomials, every
symmetric polynomial of $\gamma_1, \dots, \gamma_d$ is a
polynomial in the elementary symmetric functions of
$\gamma_1,\dots, \gamma_d$, with integers coefficients. In
particular, we can write each coefficient $\beta_j$ as a
polynomial over the integers in the variables $\alpha_1, \dots,
\alpha_d \in \R$.
\end{proof}

Since $f(x^p)=\sum_{n=0}^\infty a_{n} x^{pn}$ is again rational,
and $U_p(f(x^p))=f(x)$, we observe that the map $U_p$ is
surjective. On the other hand, for $1\le j<p$ the rational
function $x^j f(x^p)$ is in the kernel of $U_p$, so $U_p$ is not
injective. Thus, the kernel of $U_p$ is clearly infinite
dimensional. Although there is no left inverse, there is a right
inverse for $U_p$ given by the map $f(x)\mapsto f(x^p)$.

It is trivial to check that our Hecke operators form a commutative
algebra. In particular, $U_n = U_{p_1^{\alpha_1}}\circ \cdots
\circ U_{p_m^{\alpha_m}}$ whenever $n=p_1^{\alpha_1}\cdots
p_m^{\alpha_m}$.

%%%%%%%%%%%%%%%%%%%%%%%%%%%%%%%%%%%%%%%%%%%%%%%%%%%%%%%%%%%%%%%%%%%%
\section{The point spectrum of $U_p$}\label{sec_Spectrum}
%%%%%%%%%%%%%%%%%%%%%%%%%%%%%%%%%%%%%%%%%%%%%%%%%%%%%%%%%%%%%%%%%%%%

We begin by showing that the Hecke operator $U_p$ almost commutes
with the operator $x\partial_x$, up to a factor of $p$.  This
result becomes useful because it allows us to easily construct
infinitely many eigenfunctions from each known eigenfunction, by
iteration of the operator $x\partial_x$.

\begin{lemma}
For every positive integer $p$ we have
\begin{equation*}
U_p (x\partial_x) - p\,(x\partial_x)U_p = 0.
\end{equation*}
\end{lemma}
\begin{proof}
Let $f = \sum_{n=0}^\infty a_n x^n$. Then
\begin{gather*}
(x\partial_x) f = \sum_{n=1}^\infty n a_n x^n, \intertext{and it
follows that} U_p[(x\partial_x) f] = \sum_{n=1}^\infty(pn) a_{pn}
x^n  = p (x\partial_x) U_p f.
\end{gather*}
\end{proof}

%%%%%%%%%%%%%%%%%%%%%%%%%%%%%%%%%%%%%%%%%%%%%%%%%%%%%%%%%%%%%%%%%%%%
\begin{lemma}\label{IteratedEF}
Let $\lambda\not=0$ and $k\in\N$. If $U_p f = \lambda f$, then
\[ U_p[(x\partial_x)^k f]= (p^k\lambda)(x\partial_x)^k f.\]
In other words, if $\lambda$ is an eigenvalue of $U_p$, then so is
$p^k\lambda$ for every $k\in\N$, with the corresponding
eigenfunction $(x\partial_x)^k f$.
\end{lemma}
\begin{proof}
We proceed by induction on $k$.
\begin{align*}
U_p[(x\partial_x) f]
 &= p (x\partial_x) U_p f = p\lambda (x\partial_x)f
\intertext{Assuming the statement for $k-1$, we get}
U_p[(x\partial_x)^k f] &= U_p [(x\partial_x)(x\partial_x)^{k-1} f]
= p\,(x\partial_x) U_p [(x\partial_x)^{k-1} f]\\
&= p\,(x\partial_x)[p^{k-1}\lambda (x\partial_x)^{k-1}f] \\
&= (p^k \lambda) (x\partial_x)^k f.
\end{align*}
\end{proof}

%%%%%%%%%%%%%%%%%%%%%%%%%%%%%%%%%%%%%%%%%%%%%%%%%%%%%%%%%%%%%%%%%%%%

We now give an important family of eigenfunctions that have the
eigenvalues $\pm 1$ for each $U_p$. The function
$f(x)=\frac1{1-x}$ trivially satisfies $U_p f = f$ for every
positive integer $p$. However, it is less trivial to find
eigenfunctions for the eigenvalue $\lambda=-1$.

\begin{example}
For every integer $p>1$ we explicitly give an eigenfunction $f_p$
satisfying $U_p f_p = -f_p$. If $p$ is even, consider
\begin{equation}\label{evenEF}
 f_p(x)= \frac{x-x^p}{1-x^{p+1}} =
  \sum_{n=0}^\infty x^{(p+1)n+1} - \sum_{n=0}^\infty x^{(p+1)n+p}.
\end{equation}
Then, using the change of variables $j=\frac{(p+1)n+1}{p}$ and
$k=\frac{(p+1)n}{p}$ we get
\begin{align*}
U_p f_p(x) &= \sum_{\substack{pn+n+1\equiv 0 \;(p)}}
x^{\frac{(p+1)n+1}{p}}
  -\sum_{\substack{pn+n\equiv 0 \;(p)}} x^{\frac{(p+1)n}{p}+1} \\
&= \sum_{\substack{pj\equiv 1 \;(p+1)}} x^{j}
  -\sum_{\substack{pk\equiv 0 \;(p+1)}} x^{k+1} \\
&= \sum_{\substack{j\equiv p \;(p+1)}} x^{j}
  -\sum_{\substack{k\equiv 0 \;(p+1)}} x^{k+1} \\
&= \sum_{m=0}^\infty x^{(p+1)m+p}
  - \sum_{m=0}^\infty x^{(p+1)m+1} = - f_p(x).
\end{align*}

If $p$ is odd, then we can write $p-1= q\ell$ with integers $q$ and
$\ell$ such that $q$ is even and $\ell$ is odd.
In this case, the function
\begin{equation}\label{oddEF}
 f_{p}(x) =\frac{x}{1+x^q} = \sum_{n=0}^\infty (-1)^n x^{qn+1}
\end{equation}
satisfies $U_{p} f_p=-f_p$ for $p=q\ell+1$. In fact,
\begin{align*}
U_{p}\left(\frac{x}{1+x^{q}}\right)
 &= \sum_{\substack{qn+1\equiv 0 \;(p)}}
    (-1)^n x^{(qn+1)/p}\\
 &= \sum_{\substack{pk\equiv 1 \;(q)}}
    (-1)^{(pk-1)/q} x^{k}.
\intertext{Since $p\equiv 1 \bmod{q}$ it follows that $k\equiv 1
\bmod{q}$, so we write $k=qn+1$ and get}
U_{p}\left(\frac{x}{1+x^{q}}\right)
 &= \sum_{n=0}^\infty (-1)^{pn+\ell} x^{qn+1} \\
 &= -\sum_{n=0}^\infty (-1)^{n} x^{qn+1}
  = -\frac{x}{1+x^q}
\end{align*}
since $p$ and $\ell$ are both odd.
\end{example}

%%%%%%%%%%%%%%%%%%%%%%%%%%%%%%%%%%%%%%%%%%%%%%%%%%%%%%%%%%%%%%%%%%%%
The following lemma gives a first glimpse into the spectrum of
$U_p$.

\begin{lemma}\label{IntegerSpectrum}
For every integer $p\ge 2$, we have
\[ \{\pm p^{k} \st k\in\N\}\cup \{0\} \subset\spec(U_p). \]
\end{lemma}
\begin{proof}
Since $x/(1-x^p)$ is in the kernel of $U_p$, it follows that $0$
belongs to $\spec(U_p)$. The previous example gives eigenfunctions
of $U_p$ for $\lambda=1$ and $\lambda=-1$. Invoking
Lemma~\ref{IteratedEF} with these eigenfunctions, the iterated
operator $(x\partial_x)^k$ provides us with the eigenvalues
$\lambda = \pm p^k$.
\end{proof}

%%%%%%%%%%%%%%%%%%%%%%%%%%%%%%%%%%%%%%%%%%%%%%%%%%%%%%%%%%%%%%%%%%%%
Concerning the eigenfunctions of $U_p$ we have the following basic
equivalence condition.

\begin{lemma}\label{f-splitting}
Let $\lambda\not=0$.  $U_p f = \lambda f$ if and only if
\begin{equation}\label{f-identity}
 \lambda f(x^p) = \frac1p \sum_{j=0}^{p-1} f(e^{\frac{2\pi i}{p}j} x).
\end{equation}
\end{lemma}
\begin{proof}
Let $f= \sum_{n=0}^\infty a_{n} x^{n}$ be an eigenfunction associated
to $\lambda$. Then,
\begin{align*}
\lambda f(x^p) = (U_pf)(x^p) &= \sum_{k=0}^\infty a_{pk} x^{pk}\\
&= \frac1p \sum_{n=0}^\infty a_n
   (\sum_{j=0}^{p-1} e^{\frac{2\pi i}{p}nj}) x^n \\
&= \frac1p \sum_{j=0}^{p-1}
   \sum_{n=0}^\infty a_n (e^{\frac{2\pi i}{p}j} x)^n \\
&= \frac1p \sum_{j=0}^{p-1} f(e^{\frac{2\pi i}{p}j} x).
\end{align*}
Here we have used
\[ \sum_{j=0}^{p-1} e^{\frac{2\pi i}{p}nj} =
  \begin{cases}
  p& \text{for } n\equiv 0 \pmod{p}\\
  0& \text{otherwise }
  \end{cases}. \]
Now let \eqref{f-identity} be satisfied. Applying $U_p$ to both
sides of the equation we get
\begin{equation*}
 \lambda U_p(f(x^p))=\frac1p
 \sum_{j=0}^{p-1} U_p(f(e^{\frac{2\pi i}{p}j}x)) =U_p f(x)
\end{equation*}
since, for every $j$,
\[ U_p(f(e^{\frac{2\pi i}{p}j}x))=U_p\left(\sum_{n=0}^\infty a_n
  (e^{\frac{2\pi i}{p}nj}) x^n\right)= \sum_{n=0}^\infty a_{pn} x^n
  = U_p f(x). \]
Using the identity $U_p(f(x^p))=f(x)$, we conclude that $U_p
f=\lambda f$.
\end{proof}

%%%%%%%%%%%%%%%%%%%%%%%%%%%%%%%%%%%%%%%%%%%%%%%%%%%%%%%%%%%%%%%%%%%%
\begin{lemma}\label{GivenDenominator}
Let $f(x)=A(x)/B(x)$ be an eigenfunction of $U_p$. Then
\begin{equation}\label{B-identity}
 B(x^p)=\prod_{j=0}^{p-1} B(e^{\frac{2\pi i}{p}j} x).
\end{equation}
\end{lemma}

\begin{proof}
We use the identity \eqref{f-identity}
\begin{equation*}
 p\lambda\frac{A(x^p)}{B(x^p)} = p\lambda f(x^p)
 = \sum_{j=0}^{p-1} f(e^{\frac{2\pi i}{p}j} x) = \sum_{j=0}^{p-1}
 \frac{A(e^{\frac{2\pi i}{p}j} x)}{B(e^{\frac{2\pi i}{p}j} x)}
\end{equation*}
and compare the denominators.
\end{proof}

%%%%%%%%%%%%%%%%%%%%%%%%%%%%%%%%%%%%%%%%%%%%%%%%%%%%%%%%%%%%%%%%%%%%
\begin{theorem}\label{MatrixB}
If $\lambda$ is an eigenvalue of $U_p$ and $f(x)=A(x)/B(x)$ is a
corresponding eigenfunction, then $\lambda$ is an eigenvalue of a
$d\times d$ matrix $\mathfrak{B}$ determined by the coefficients
of the polynomial $B(x)$.
\end{theorem}
\begin{proof}
Using \eqref{f-identity} and \eqref{B-identity} we get the identity
\begin{equation*}
 \lambda A(x^p) = \frac1p \sum_{j=0}^{p-1}
 \bigg(\prod_{\ell\not= j} B(\zeta_p^\ell x) \bigg) A(\zeta_p^j x),
\end{equation*}
where $\zeta_p = e^{\frac{2\pi i}{p}}$. Writing
$A(x)=\sum_{k=0}^{d-1} c_k x^k$,  the identity above becomes
\begin{align*}
 \lambda\sum_{k=0}^{d-1} c_k x^{pk} &= \frac1p \sum_{j=0}^{p-1}
 \bigg(\prod_{\ell\not= j} B(\zeta_p^\ell x) \bigg)
  \sum_{k=0}^{d-1} c_k (\zeta_p^j x)^{k} \\
 &= \sum_{k=0}^{d-1} \bigg(\sum_{m} \beta_{k,m} c_m \bigg) x^{pk}
\end{align*}
by rearranging the sum on the right-hand side. Comparing the coefficients
gives us
\[ \lambda c_k = \sum_{m=0}^{d-1} \beta_{k,m} c_m \]
for every $k=0,\dots,d-1$.
Finally, we conclude that $\lambda$ is an eigenvalue of the
matrix $\mathfrak{B} = (\beta_{k,m})_{k,m=0,\dots,d-1}$. We observe that
the $\beta_{k,m}$ are complex polynomials in the coefficients of the
denominator $B(x)$.
\end{proof}

%%%%%%%%%%%%%%%%%%%%%%%%%%%%%%%%%%%%%%%%%%%%%%%%%%%%%%%%%%%%%%%%%%%%

This theorem is very useful because it allows us to explicitly construct
eigenfunctions with a given denominator by computing eigenvectors of a
\emph{finite} matrix.  As a direct consequence of this theorem, we obtain a
converse to Lemma~\ref{GivenDenominator}.
\begin{corollary}
Given any integer $p$ and a denominator $B(x)$ satisfying the
identity \eqref{B-identity}, there is a numerator $A(x)$ such that
$f(x)=A(x)/B(x)$ is an eigenfunction of $U_p$.
\end{corollary}

%%%%%%%%%%%%%%%%%%%%%%%%%%%%%%%%%%%%%%%%%%%%%%%%%%%%%%%%%%%%%%%%%%%%

\begin{theorem}[Involution Property]\label{B-structure}
If $f(x)=A(x)/B(x)$ is an eigenfunction of $U_p$, then all roots
of $B(x)$ are roots of unity and we have the identity
\begin{equation*}
 x^d B(\tfrac{1}{x}) = (-1)^d B(x).
\end{equation*}
Moreover, if $U_p f = \lambda f$ with $\lambda \not= 1$, then
\[ U_p (f(\tfrac{1}{x})) = \lambda f(\tfrac{1}{x}).  \]
That is,  $f(\tfrac{1}{x})$ is another eigenfunction of $U_p$ with
the same eigenvalue $\lambda$, and with the same denominator
$B(x)$.
\end{theorem}
\begin{proof}
Write
\begin{equation*}
B(x)=1+\alpha_1 x +\cdots +\alpha_d x^d=\prod_{k=1}^d (1-\gamma_k x).
\end{equation*}
The identity \eqref{B-identity} yields
\begin{align*}
\prod_{k=1}^d (1-\gamma_k x^p)
&=\prod_{j=0}^{p-1} \prod_{k=1}^d (1-\gamma_k(e^{\frac{2\pi i}{p}j} x)) \\
&=\prod_{k=1}^d \prod_{j=0}^{p-1} (1-\gamma_k(e^{\frac{2\pi i}{p}j} x)) \\
&=\prod_{k=1}^d (1-\gamma_k^p x^p)
\end{align*}
which implies
$\{\gamma_1,\dots,\gamma_d\}=\{\gamma_1^p,\dots,\gamma_d^p\}$.
So the second set is a permutation of the first set. This
permutation breaks up into a disjoint product of cycles.
Consider now a fixed cycle in this decomposition, say of
length $\ell$. By iterating through the cycle, we easily see
that each of the roots in this cycle must satisfy the equation
\begin{equation}\label{RootOfUnity}
x^{p^\ell -1} = 1.
\end{equation}
Thus they are all roots of unity (different from $-1$ if $p$ is
even). Since we can do this for each cycle, all of the
$\gamma_j$'s are in fact roots of unity.

As a consequence, we get $\prod_{k=1}^d \gamma_k=1$ and
$\gamma_k^{-1} = \bar\gamma_k$. Therefore,
\begin{align*}
x^d B(\tfrac{1}{x}) &= x^d \prod_{k=1}^d (1- \gamma_k/x)
= (-1)^d \prod_{k=1}^d (\gamma_k-x) \\
&= (-1)^d \prod_{k=1}^d \gamma_k (1 - \gamma_k^{-1}x) \\
&= (-1)^d \prod_{k=1}^d (1 - \gamma_k^{-1}x)
 = (-1)^d \prod_{k=1}^d (1 - \bar\gamma_k x) \\
&= (-1)^d \prod_{k=1}^d (1 - \gamma_k x)
 = (-1)^d B(x).
\end{align*}

To prove the last claim, let $f(x)=A(x)/B(x)$ be an eigenfunction
of $U_p$ with eigenvalue $\lambda\not=1$. We will verify that the
function $g(x) = f(\frac{1}{x})$ satisfies the condition
\eqref{f-identity}. The identity $x^d B(\tfrac{1}{x})=(-1)^d B(x)$
gives
\begin{equation*}
 g(x) = \frac{A(\frac1x)}{B(\frac1x)} =
 (-1)^d\; \frac{x^d A(\frac1x)}{B(x)},
\end{equation*}
where $x^d A(\frac1x)$ is a polynomial of degree less than $d$ since
$\lambda\not=1$ implies $A(0)=0$. Thus $g$ belongs to $\rational$
and has the same denominator as $f$. Now,
\begin{align*}
\sum_{j=0}^{p-1} g(\zeta_p^j x)
  = \sum_{j=0}^{p-1} \frac{A\Big(\frac{1}{\zeta_p^j x}\Big)}
  {B\Big(\frac{1}{\zeta_p^j x}\Big)}
  = \sum_{k=0}^{p-1} \frac{A(\zeta_p^{k}y)}{B(\zeta_p^{k}y)}
  = \sum_{k=0}^{p-1} f(\zeta_p^{k}y),
\end{align*}
where $\zeta_p = e^{\frac{2\pi i}{p}}$ and $y = \frac{1}{x}$.  In the
latter equality we used the fact that $1/\zeta_p^j = \zeta_p^{k}$ with $k=p-j$.
Using the condition \eqref{f-identity} for $f$, we obtain
\begin{equation*}
 \frac1p \sum_{j=0}^{p-1} g(\zeta_p^j x) =
 \frac1p \sum_{k=0}^{p-1} f(\zeta_p^{k}y) = \lambda f(y^p) =
 \lambda f(\tfrac{1}{x^p}) = \lambda g(x^p).
\end{equation*}
Thus $g(x)$ satisfies \eqref{f-identity} and $f(\frac{1}{x})$ is
therefore another eigenfunction of $U_p$ with the same eigenvalue
$\lambda$.
\end{proof}

It is interesting to note the analogy here with the classical
Fricke involutions on Riemann surfaces. The preceding theorem
gives us an involution on the vector space of eigenfunctions, and
supplies us with an easy construction of new eigenfunctions from
known ones.

\begin{corollary}\label{Cor4-2}
Let $f$ be any rational function such that $U_p f = \lambda_p f$
for some integer $p$, and let $L$ be the level of $f$. Then $p$ is
relatively prime to $L$.
\end{corollary}
\begin{proof}
The relation \eqref{RootOfUnity} appearing in the proof of 
previous theorem tells us that $p^l -1$ must be a 
multiple of the level $L$.  Therefore $p(p^{l-1}) -1 = mL$ and
we see that $p$ is relatively prime to $L$. 
\end{proof}

\begin{definition}
Fix $p$ and suppose we have an eigenfunction $f$ for $U_p$ with
eigenvalue $\lambda_p$. If $\lambda_p \not= 0$, we define
\begin{equation*}
 \chi_f(p) = \frac{\lambda_p}{|\lambda_p|}.
\end{equation*}
If $\lambda_p = 0$, then we define $\chi_f(p) = 0$.
\end{definition}

In Section~\ref{sec_Simultaneous} we will prove that, if $f$ is a
simultaneous eigenfunction for an appropriate family $\{U_p \st
p\in S\subset\N\}$, then $\chi_f$ is a \emph{character}.

\begin{lemma}\label{RealSpectrum}
For every integer $p$, $\spec(U_p) \subset \R$. In particular, for
every $f$ with $U_p f=\lambda f$, $\lambda\not=0$, we have
$\chi_f(p)=\pm 1$.
\end{lemma}
\begin{proof}
Let $f$ be an eigenfunction of $U_p$ for $\lambda\not=0$. Since
$f(\bar\rho x) = \overline{f(\rho x)}$ for any complex number
$\rho$, we have
\[ f(e^{\frac{2\pi i}{p}(p-j)} x)= f(e^{-\frac{2\pi i}{p}j} x)
  = \overline{f(e^{\frac{2\pi i}{p}j} x)} \]
for every $j=1,\dots,p-1$.
Thus $f(e^{\frac{2\pi i}{p}j} x) + f(e^{\frac{2\pi i}{p}(p-j)} x)$
is always real, which implies that the right-hand side of
\eqref{f-identity} is real. Consequently, $\lambda$ is real as well.
\end{proof}

%%%%%%%%%%%%%%%%%%%%%%%%%%%%%%%%%%%%%%%%%%%%%%%%%%%%%%%%%%%%%%%%%%%%
The next theorem gives us the complete structure of the spectrum
for every Hecke operator $U_p$. We note that the proof of the main
structure theorem for eigenfunctions,
Theorem~\ref{EFunctionStructure} in the next section, bootstraps
the proof of the following theorem.

\begin{theorem}[The Spectrum]\label{EValueStructure}
Let $p$ be any positive integer greater than $1$. Then the point
spectrum of $U_p$ on the vector space $\rational$ is
\[ \spec(U_p) = \{\pm p^{k} \st k\in\N\}\cup \{0\}. \]
\end{theorem}
\begin{proof}
Let $f(x)=\sum a_n x^n$ satisfy $U_pf=\lambda f$, $\lambda\not=0$.
Thus $a_{pn} = \lambda a_n$  and by iteration
\begin{equation}\label{Eq1}
\,a_{p^k n}=\lambda^k a_n
\end{equation}
for every $k\in\N$. Let $B(x)= \prod_{j=1}^d (1-\gamma_j x)$. From
Theorem~\ref{B-structure} we know that the $\gamma_j$'s must be
roots of unity.  We now compare the asymptotics of the
coefficients from their closed form as a sum of polynomials in $n$
times roots of unity. Namely, from Theorem~\ref{RationalFunctions}
we know that
\begin{equation*}
a_n = \sum_{j=1}^\ell  c_j n^{m_j-1} \gamma_j^n,
\end{equation*}
where each $m_j$ is the multiplicity of the root $\gamma_j$. For
the purposes of using asymptotics, we let $\kappa$ denote the
largest exponent $m_j$ in this representation of $a_n$.  We
collect together all of the terms that correspond to this largest
exponent $\kappa$, and label the remaining terms by $R(n)$. Thus
we may write
\begin{equation}\label{Eq2}
 a_n  = n^{\kappa-1} \sum_{j=1}^{\ell_1} C_{j}\gamma_{\sigma(j)}^{n}+
 R(n),
\end{equation}
for some constants $C_j$ and some permutation $\sigma$. It follows
from Equation~\ref{Eq1} above (using the assumption that $\lambda
\not=0$) that
\begin{equation}\label{Eq3}
 a_n = \frac{a_{p^k n}}{\lambda^k} =
 \left(\frac{p^{\kappa-1}}{\lambda}\right)^k
 n^{\kappa-1} \sum_{j=1}^{\ell_1} C_{j}\gamma_{\sigma(j)}^{p^k n}+
 \frac{R(p^k n) }{\lambda^k}.
\end{equation}

We first claim that $|\frac{p^{\kappa-1}}{\lambda}| \leq 1$. To see
this note that all of the terms in $ \frac{R(p^k n)}{\lambda^k}$
contain exponential terms in $k$ of the form
$(\frac{p^{m_j-1}}{\lambda})^k$, which are strictly of smaller growth
than $(\frac{p^{\kappa-1}}{\lambda})^k$.  Since the left-hand side of
the equation above is $a_n$, and in particular independent of $k$, we
cannot have $|\frac{p^{\kappa-1}}{\lambda}| > 1$, unless perhaps the
leading sum vanishes for all $k$.  We now argue that there is a
subsequence of $k$'s for which the leading sum $\sum_{j=1}^{\ell_1}
C_{j}\gamma_{\sigma(j)}^{p^k n}$ does not vanish.

To this end note that, by Corollary~\ref{Cor4-2}, $p$ is relatively
prime to $L$, where $L$ is the least common multiple of the orders (as
roots of unity) of all the poles of $f$.  
Using $p^{\phi(L)} \equiv 1 \pmod L$,
we now let the index $k$ approach infinity through the
subsequence $k' = m \phi(L) $, where $m \in \N $.
Thus $\gamma_{\sigma(j)}^{p^{k'}} = \gamma_{\sigma(j)}$, for all 
$1\leq j \leq \ell_1$, and we have
\begin{equation}\label{Eq4}
 a_n =
 \left(\frac{p^{\kappa-1}}{\lambda}\right)^{k'}
 n^{\kappa-1} \sum_{j=1}^{\ell_1} C_{j}\gamma_{\sigma(j)}^{n}+
 \frac{R(p^{k'} n) }{\lambda^{k'}}.
\end{equation}
Thus the term with the largest exponent cannot vanish as
$k \to \infty$ through the given subsequence, so that have shown that
$|\frac{p^{\kappa-1}}{\lambda}| \leq 1$.

On the other hand, we cannot have $|\frac{p^{\kappa-1}}{\lambda}| < 1$, for
then all of the terms on the right hand side of Equation~\ref{Eq3}
would tend to $0$ as $k \rightarrow \infty$, contradicting $a_n \not= 0$.
Therefore
$|\frac{p^{\kappa-1}}{\lambda}| = 1$ and thus $\lambda=\pm
p^{\kappa-1}$ since $\lambda$ is real by Lemma~\ref{RealSpectrum}.
Together with the inclusion from Lemma~\ref{IntegerSpectrum} we
finally get the assertion.
\end{proof}

%%%%%%%%%%%%%%%%%%%%%%%%%%%%%%%%%%%%%%%%%%%%%%%%%%%%%%%%%%%%%%%%%%%%
\section{A structure theorem for eigenfunctions}\label{sec_Structure}
%%%%%%%%%%%%%%%%%%%%%%%%%%%%%%%%%%%%%%%%%%%%%%%%%%%%%%%%%%%%%%%%%%%%

By refining the proof of Theorem~\ref{EValueStructure} further, we
can get a very useful structure theorem for eigenfunctions.  We can
subsequently draw several interesting conclusions that resemble
ideas from automorphic forms.  In particular, we will define a
\emph{weight} and a \emph{level} for eigenfunctions, and show in
Section~\ref{sec_PGR} that we have a finite dimensional eigenspace
for a fixed level $L$ \emph{and} weight $\kappa$.

\begin{definition}
Given an eigenfunction $f$ of $U_p$, we know by
Theorem~\ref{B-structure} that its poles are all roots of unity.
We define the \emph{level} $L$ of $f$ as the least common multiple
of the orders of all these roots of unity; thus each pole
$\gamma_j$ is a root of unity $e^{\frac{2\pi i \ell_j}{L}}$ for
some integer $\ell_j$.
\end{definition}

We observe that when the level of $f$ is $L$, the smallest group
containing all of the poles of $f$ is simply the group $\mu_L$ of
$L$'th roots of unity.  We now give a structure theorem that
simplifies the analysis of eigenfunctions.

\begin{theorem}[Structure Theorem] \label{EFunctionStructure}
Let $f(x)=A(x)/B(x)=\sum a_n x^n$ be an eigenfunction of $U_p$ for
some integer $p>1$, associated to an eigenvalue $\lambda_p
\not=0$. If $B(x)= \prod_{j=1}^d (1-\gamma_j x)$, then there is an
integer $\kappa$ dividing the degree $d$, and an integer $L$ such that
\begin{equation*}
 a_n = n^{\kappa-1}\sum_{j=1}^{d/\kappa} C_{j}\,
 e^{\frac{2\pi i  \ell_j}{L}n},
\end{equation*}
where each pole of $f$ is given by $\gamma_j =e^{\frac{2\pi i
\ell_j}{L}}$, $\ell_j\in\N$, and the constants $C_{j}\in \C$ are determined
by the initial conditions of the linear recurrence sequence
$\{ a_n \}$. We note that each pole $\gamma_j$ must occur with the
same multiplicity $\kappa$.
\end{theorem}
\begin{proof}
We begin with the identity \eqref{Eq.a_n} derived in the proof of
Theorem~\ref{EValueStructure}, namely
\begin{align*}
 a_n &=\left(\frac{p^{\kappa-1}}{\lambda}\right)^k
 n^{\kappa-1} \sum_{j=1}^\ell C_{j}\gamma_{\sigma(j)}^{p^k n}+ R(k)\\
 &=\chi^k_f(p)\, n^{\kappa-1} \sum_{j=1}^{\ell} C_{j}\,
 e^{\frac{2\pi i  \ell_j}{L}p^k n}+ R(k) \\
 &=(\pm 1)^k n^{\kappa-1} \sum_{j=1}^{\ell} C_{j}\,
 e^{\frac{2\pi i  \ell_j}{L}p^k n}+ R(k),
\end{align*}
where we used Lemma~\ref{RealSpectrum} to express $\chi_f(p)=\pm 1$.
We now claim that $R(k)$ is identically zero. To see this, note that the
sum on the right-hand side is a periodic function in $k$ which implies
that $R(k)$ is also periodic in $k$.

On the other hand, recall that all of the terms in $R(k)$ contain exponential
terms of the form $(\frac{p^{m_j-1}}{\lambda})^k$ with $m_j<\kappa$.
Since $p^{m_j-1}<p^{\kappa-1}$, we have
$|\frac{p^{m_j-1}}{\lambda}| < |\frac{p^{\kappa-1}}{\lambda}|=1$, so that
$|\frac{p^{m_j-1}}{\lambda}|\to 0$.
Thus $R(k)\to 0$ as $k\to\infty$ and we conclude from the periodicity of
$R(k)$ that it is the zero function.
\end{proof}

\begin{definition}
We call the $\kappa$ appearing in the previous theorem the {\em
weight} of the eigenfunction $f$. Note that $\kappa$ is
independent of $p$. That is, if $f$ is an eigenfunction of any
other $U_q$, then $U_q f = \pm q^{\kappa-1} f$.
\end{definition}

The motivation for this terminology comes from the weight of the classical
Eisenstein series in automorphic forms.  With hindsight, we call
the operator $x\partial_x$ the weight-raising operator, since it
takes a weight $\kappa$ eigenfunction to a weight $\kappa + 1 $ eigenfunction.

\begin{remark}\label{DegreeDenominator}
There is a precise connection between the weight $\kappa$, the level
$L$ and the degree of the denominator of an eigenfunction
$f(x)=A(x)/B(x)$.  Namely, $\deg(B(x))\in \{\kappa, 2\kappa, \dots,
L\kappa\}$.
\end{remark}

At this stage we note that, given $p>1$, every eigenvalue of $U_p$
must have the form
\begin{equation}\label{Eigenvalues}
 \lambda =  \chi_f(p) \  p^{\kappa -1},
\end{equation}
where $\chi_f(p)=\pm 1$. In the corollaries that follow we will show that
$\chi_f$ has the properties of a multiplicative character.
%%%%%%%%%%%%%%%%%%%%%%%%%%%%%%%%%%%%%%%%%%%%%%%%%%%%%%%%%%%%%%%%%%%%

\begin{corollary}\label{Cor4-1}
Let $f=\sum a_n x^n$ be a rational function such that $U_p f =
\lambda_p f$ for some integer $p$, and let $L$ be the level of
$f$. Then
\begin{enumerate}
\item[(i)] If $a_0 = 0$, then $\chi_f(L) = 0$ and $U_{mL}(f) = 0$
for every positive integer $m$. That is, $a_{nL} = 0$ for every
$n$. \item[(ii)] If $a_0 \not= 0$, then $\chi_f(p) =1$, so that
$U_p(f) = f$.
\end{enumerate}
\end{corollary}
\begin{proof}
To show (i), we use the structure theorem directly, keeping in mind
that by definition of $L$ we have $\gamma_j^L =1$ for all of the roots
$\gamma_j$.
$$
a_{nL} = (nL)^{\kappa -1} \sum_j C_j \gamma_j^{nL} =
(nL)^{\kappa -1} \sum_j C_j = (nL)^{\kappa -1}a_0 = 0.
$$
Thus by definition of $U_{mL}$, we have $U_{mL}(f) = 0$.

To show (ii) we simply recall that $a_0 \lambda = a_0$, so that
$\lambda =1= \chi_f(p) p^{\kappa-1}$, and we have the required
result.
\end{proof}

In other words any eigenfunction $f$,  with eigenvalue $\lambda
\not= 1$ and  level $L$, must lie in the kernel of $U_L$ and thus
$f$ has an infinite arithmetic progression of zeros among its
Taylor coefficients.

\begin{corollary}\label{Cor4-3}
For any positive integer $n$, we have the following formulas for $\chi_f(p)$:
\begin{enumerate}
\item[(i)] $\chi_f(p) \sum_j C_j \gamma_j^n = \sum_j C_j
\gamma_j^{pn}$. \item[(ii)] $\chi_f(p+mL) = \chi_f(p)$ \text{\ for
all positive integers} $m$. \item[(iii)] Suppose that $f$ is a
simultaneous eigenfunction for all $U_p$.  Then $f$ is not in the
kernel of any operator $U_p$ if and only if $\chi_f(p) = 1$ for
all $p$.
\end{enumerate}
\end{corollary}
\begin{proof}
We begin with the identity
\[ a_{pn} = (pn)^{\kappa-1}
  \sum_{j=1}^{d/\kappa} C_{j}\, e^{\frac{2\pi i  \ell_j}{L}pn}. \]
On the other hand,
\[ a_{pn} = \chi_f(p) p^{\kappa-1} a_n =
 \chi_f(p) p^{\kappa-1}  n^{\kappa-1}\sum_{j=1}^{d/\kappa} C_{j}\,
 e^{\frac{2\pi i  \ell_j}{L}n}. \]
Equating both right-hand sides, we get the desired identity (i).

To prove (ii), we simply pick an $n$ for which $0 \not= a_n = n^{\kappa
-1}\sum_j C_j \gamma_j^n$, whence
\[ \chi_f(p + mL) =
  \frac{ \sum_j C_j \gamma_j^{(p+mL)n}}{\sum_j C_j \gamma_j^n}
 =\frac{\sum_j C_j \gamma_j^{pn}}{\sum_j C_j \gamma_j^n } = \chi_f(p). \]

To prove (iii), first assume that $f$ is not in the kernel of any
operator $U_p$. We observe that Corollary~\ref{Cor4-2} implies
$(p, L) = 1$.  Since $L$ is now relatively prime to   \textit{all}
integers $p$, we must have $L=1$ which implies $\gamma_j=1$ for
all $j$. Using $L=1$ in part (i) now gives us
\[ \chi_f(p) \sum_j C_j  = \sum_j C_j, \]
and we conclude that either $\chi_f(p)=1$ for each $p$ (and we're
done) or else $\sum_j C_j =0$.  To see that the latter case never
occurs, we observe that $\sum_j C_j =0$ means that $a_0 = 0$,
which by Corollary~\ref{Cor4-1} (i) above implies that $f$ is in
the kernel of $U_L$, a contradiction.

To prove (iii) in the other direction, let $\chi_f(p) = 1$ for all
$p$. Then $U_p f  = \chi_f(p) p^{\kappa-1} f = p^{\kappa-1} f$,
for all $p$, so that f is trivially not in the kernel of any Hecke
operator $U_p$.
\end{proof}

\begin{corollary}\label{Cor4-4}
Whenever $f$ is a simultaneous eigenfunction of two distinct
operators $U_m$ and $U_n$, then $f$ is also an eigenfunction of
$U_{mn}$ and in particular we obtain the identity
\[ \chi_f(mn) = \chi_f(m) \chi_f(n). \]
\end{corollary}
\begin{proof}
We first compute
\[ U_{mn}(f) = U_m(U_n f) = U_m( \chi_f(n)n^{\kappa-1} f ) =
  \chi_f(m) m^{\kappa-1} \chi_f(n) n^{\kappa-1} f, \]
whence $f$ is indeed an eigenfunction of $U_{mn}$. From the
structure theorem, we therefore obtain
\[ U_{mn}(f) = \chi_f(mn) (mn)^{\kappa-1} f, \]
where we note that the same weight $\kappa$ appears in both computations,
due to the fact that the weight $\kappa$ depends only on $f$ and not on
the Hecke operator.  The result follows by comparing the two
equalities above.
\end{proof}

\begin{corollary}\label{Cor4-5}
Let $f$ be any rational function such that $U_p f = \lambda_p f$
for some integer $p$, and let $L$ be the level of $f$. Then
\begin{equation*}
 U_{p+mL}(f) = \lambda_{p+mL}f
\end{equation*}
for every positive integer $m$. In addition, whenever $f$ is an
eigenfunction of a single operator $U_p$ it is also an
eigenfunction of an infinite collection of operators $U_q$ with
prime index $q$.
\end{corollary}
\begin{proof}
For any positive integer $m$, it follows from the structure theorem
with $\gamma_j = e^{\frac{2\pi i \ell_j}{L}}$ that
\begin{align*}
 a_{(p+mL)n} &=
 (p+mL)^{\kappa-1}n^{\kappa-1}
 \sum_{j=1}^{d/\kappa} C_{j}\, \gamma_j^{(p+mL)n} \\
 &= (p+mL)^{\kappa-1}n^{\kappa-1}
 \sum_{j=1}^{d/\kappa} C_{j}\, \gamma_j^{pn} \\
 &= (p+mL)^{\kappa-1}n^{\kappa-1} \chi_f(p)
 \sum_{j=1}^{d/\kappa} C_j \gamma_j^n \\
 &= (p+mL)^{\kappa-1} \chi_f(p+mL) n^{\kappa-1}
 \sum_{j=1}^{d/\kappa} C_j \gamma_j^n \\
 &=  (p+mL)^{\kappa-1}  \chi_f(p+mL) a_n,
\end{align*}
where the third equality is part (i) of Corollary~\ref{Cor4-3},
the fourth equality is part (ii) of Corollary~\ref{Cor4-3}, and
the last equality simply uses the structure theorem for
eigenfunctions. Thus, by definition, we arrive at $U_{p+mL}(f) =
\lambda_{p+mL} f$.

The second claim follows from Dirichlet's theorem on primes in
non-trivial arithmetic progressions
(see Knapp \cite[p.  189]{MR93j:11032}),
once we know that $(p, L) =1$ from Corollary~\ref{Cor4-2}.
\end{proof}

The latter result links the index $p$ of the operator with the
level $L$ of the eigenfunction in a strong way.  This connection
makes $L$ a natural candidate for grading the eigenspaces of
$U_p$, a task we take up in the following section.

We conclude this section by showing that $\chi_f$ is in fact the real
quadratic character mod $L$ if we know that the real rational function
$f$ is a simultaneous eigenfunction of sufficiently many Hecke operators.

%%%%%%%%%%%%%%%%%%%%%%%%%%%%%%%%%%%%%%%%%%%%%%%%%%%%%%%%%%%%%%%%%%%%
\begin{theorem} \label{Character}
Let $L>1$ be a given integer.
\begin{enumerate}
\item[(i)] Suppose $f(x)=\sum_{n=0}^\infty a_n x^n$ is a real
rational function of level $L$ with $a_0 = 0$, and $f(x)$ is a
simultaneous eigenfunction of the operators $U_2, U_3,  \dots,
U_L$ (i.e. $\;U_m f = \chi_f(m) m^{\kappa -1} f$ for every
$m=2,\dots,L$). Then $\chi_f$ is the real quadratic character mod
$L$, $f$ is in fact a simultaneous eigenfunction of all the Hecke
operators $U_m$, and in addition we must have
\[  f(x) = a_1 \sum_{n=0}^\infty  \chi_f(n) n^{\kappa -1} x^n. \]

It is worthwhile noting that $f$ can also be written as
\begin{equation*}
  f(x) = a_1 \ (x \partial_x)^{\kappa-1}
  \bigg(\frac{\sum_{j=1}^{L-1} \chi_f(j) x^j}{1-x^L} \bigg).
\end{equation*}
Under the same hypothesis, except with $a_0\not=0$, we conclude
that $\chi_f$ is the identity character and
$f(x)=\frac{a_0}{1-x}$.
\item[(ii)] Conversely, given any (real or
complex) character $\chi$ mod $L$, and any positive integer
$\kappa$, the rational function
\begin{equation*}
 f(x) = \sum_{n=0}^\infty \chi(n) n^{\kappa-1} x^n
\end{equation*}
satisfies $U_p f = \chi(p) p^{\kappa-1} f$ for every $p$.
\end{enumerate}
\end{theorem}

\begin{remark}
In other words, the previous theorem tells us that
when we restrict a level $L$ rational function $f$ to be a simultaneous
eigenfunction of the first $L$ Hecke operators, the function $f$ must lie
in the $1$-dimensional vector space generated by the given function.
\end{remark}

\begin{proof}
To prove (i), we start with Corollary~\ref{Cor4-3}(ii), from which
we know that the sequence of real values $\{\chi_f(1), \chi_f(2),
\chi_f(3), \dots, \chi_f(L) \}$ extends to all of $\N$ by the
periodicity of $\chi_f$ mod $L$.  We also know from
Corollary~\ref{Cor4-4} that these values are multiplicative. Thus
we have a real character mod $L$.  The fact that $f$ is a
simultaneous eigenfunction of all the Hecke operators follows from
the previous Corollary~\ref{Cor4-5}, which tells us that it
suffices to only consider those Hecke operators $U_p$ with $p$
less than or equal to the level of $f$.

Now let $f(x)=\sum_{n=1}^\infty a_n x^n$ be an eigenfunction satisfying the
hypothesis. Then $a_{p}=\chi_f(p)\, p^{\kappa-1} a_1$ for every $p=2,\dots,L$.
Thus for $n>L$, we have $a_n=a_{p+jL}=\chi_f(p+jL)(p+jL)^{\kappa-1}a_1$ by
Corollary~\ref{Cor4-5}, and so
\begin{align*}
 f(x) &= a_1 \sum_{n=1}^\infty \chi_f(n)\, n^{\kappa-1} x^n.
\end{align*}

In the case that $a_0 \not=0$, we use the fact that $\lambda_p a_0 = a_0$ for
each $2 \leq p \leq L$ to get $1 = \lambda_p = \chi_f(p) p^{\kappa -1}$. Thus
$\kappa = 1$ and $\chi_f(p) = 1$ for each such $p$. By the periodicity
of $\chi_f$, we obtain $\chi_f(n) = 1$ for all positive integers $n$, and
therefore
\begin{align*}
 f(x) &= a_0 + \sum_{n=1}^\infty \chi_f(n)a_1 x^n \\
 &= a_0 + a_1\sum_{n=1}^\infty x^n
 = \frac{a_0 + (a_1-a_0) x}{1-x} = \frac{a_0}{1-x},
\end{align*}
since $a_1=a_0$, a conclusion that follows from the fact that the
degree of the numerator must be smaller than the degree of the
denominator ($f\in\rational$).

To prove (ii), we just compute
\begin{align*}
 U_pf(x) &= \sum_{n=1}^\infty \chi(p n) (pn)^{\kappa-1} x^n \\
 &= \chi(p) p^{\kappa-1} \sum_{n=1}^\infty \chi(n) n^{\kappa-1} x^n \\
 &= \chi(p) p^{\kappa-1} f(x).
\end{align*}
Finally, we apply the weight-raising operator $(x \partial_x)$
to get
\[f(x)= (x \partial_x)^{\kappa-1}
\big(\sum_{n=1}^{\infty} \chi(n) x^n \big).  \]
Thus the second representation of $f(x)$ follows from the identity
\begin{align*}
 \sum_{n=1}^{\infty} \chi(n) x^n
 &= \sum_{j=1}^{L-1}\sum_{m=0}^{\infty} \chi(j+mL) x^{j+mL}\\
 &= \sum_{j=1}^{L-1} \chi(j) x^j
    \sum_{m=0}^{\infty} x^{mL}\\
 &= \frac{\sum_{j=1}^{L-1} \chi(j) x^j}{1-x^L},
\end{align*}
where we used the property $\chi(j+mL) =\chi(j)$.
\end{proof}

\begin{example}\label{Level7}
We give an example of some rational functions of level $7$ that
are simultaneous eigenfunctions of exactly two Hecke operators,
but not of all of them. Due to the periodicity property
$U_{p+mL}(f) = \lambda_{p+mL} f$ of Corollary~\ref{Cor4-5} it
suffices to consider only Hecke operators $U_p$ with $2\le p\le
7$.  This illustrates the interesting fact that although the
values $\chi_f$ are multiplicative, they can be restricted away
from actually being the full character mod $L$, and we thus get
``partial'' characters.

Let $f(x) = \frac{x+x^2 +x^4}{1-x^7}$, and let
$g(x) = \frac{x^3+x^5+x^6}{1-x^7}$.
It is clear that both $f$ and $g$ are of weight $1$ and level $7$, since
the poles of each function are in fact all the distinct $7$'th roots of
unity.  Furthermore, we have
\begin{align*}
 f(x) &= x+x^2+x^4 + x^8+x^9+x^{11} +x^{15}+x^{16}+x^{18}+ \cdots
\intertext{and}
 g(x) &= x^3+x^5+x^6 + x^{10}+x^{12}+x^{13} +x^{17}+x^{19}+x^{20}
  + \cdots
\end{align*}
 From the power series, it is trivial to check that $U_2 f = f$,
$U_4 f = f$, $U_2 g = g$, and $U_4 g = g$.  Thus in this example
$\chi_f(2)= \chi_f(4) = 1$ and $\chi_g(2)= \chi_g(4)=1$, but
$\chi_f(p)$ is not even defined for other values of $p$ mod $7$.
For all other $p$ mod $7$, it is easy to see from their Taylor
series that the functions $f$ and $g$ are not eigenfunctions of
any other $U_p$. We note that the numerators of eigenfunctions are
in general non-trivial polynomials and it is an interesting (and
in general difficult) problem to compute them. Furthermore this
example illustrates the involution property from
Theorem~\ref{B-structure}. Indeed, $f(\frac{1}{x})= -g(x)$.
\end{example}

%%%%%%%%%%%%%%%%%%%%%%%%%%%%%%%%%%%%%%%%%%%%%%%%%%%%%%%%%%%%%%%%%%%%
\section{A decomposition into finite dimensional eigenspaces}
\label{sec_PGR}
%%%%%%%%%%%%%%%%%%%%%%%%%%%%%%%%%%%%%%%%%%%%%%%%%%%%%%%%%%%%%%%%%%%%

Let $f(x)=\sum_{n=0}^{\infty}a_n x^n\in \rational_{qp}$ be any
real rational function in $\rational$ whose poles are roots of unity.
In this section we first show that $f$ lies in the real span of some very
simple Hecke
eigenfunctions, of the same level $L$.  We then define some finite
dimensional vector spaces of eigenfunctions that have fixed weight and
level, again by analogy with automorphic forms.

\begin{theorem} \label{PGRBasis}
Let $f(x) = \sum_{n=0}^\infty a_n x^n$ be any rational function in
$\rational_{qp}$.  Then $f$ lies in the real span of the following
eigenfunctions of $U_{L+1}$, each eigenfunction having level $L$:
\[
\left( x\partial_x  \right)^{k}  \left( \frac{x^{j} }{ 1-x^L} \right),
\]
for all non-negative integers \ $k$ and  $ 0 \leq j <L$.
\end{theorem}

\begin{proof}
 From the standard Theorem~\ref{RationalFunctions} we know that $a_n$ is
a quasi-polynomial in $n$.  That is, we have $a_n = \sum_{j=1}^L
n^{m_j} b_j(n)$, where the $b_j$'s are periodic, real-valued functions
of $n$.  Let $L$ be the least common multiple of all the periods of
the $b_j$'s.  We first observe that for any infinite
periodic sequence $b(n)$ of period $L$, we have
\begin{align*}
\sum_{n=0}^\infty b(n) x^n
 &= \frac{b(0)}{1-x^L} + \frac{b(1)x}{1-x^L}+ \dots +
 \frac{b(L-1)x^{L-1}}{1-x^L} \\
 &= \frac{ b(0) + b(1)x +\cdots + b(L-1) x^{L-1} }{ 1-x^L}.
\end{align*}
Thus, for each index $j$ we obtain
\[
 \sum_{n=0}^\infty n^{m_j} b_j(n) x^n = \left( x\partial_x  \right)^{m_j}
 \left(
 \frac{ b_j(0) + \dots + b_j(L-1) x^{L-1} }{ 1-x^L}
 \right),
\]
and consequently
\begin{align*}
 f(x) &=
 \sum_{n=0}^\infty a_n x^n = 
 \sum_{n=0}^\infty \sum_{j=1}^L n^{m_j} b_j(n) x^n \\
 &= \sum_{j=1}^L \left( x\partial_x  \right)^{m_j}
 \left(
 \frac{ b_j(0) + \dots + b_j(L-1) x^{L-1} }{ 1-x^L}
 \right)   \\
 &=  \sum_{j=1}^L \left(
  b_j(0) \left( x\partial_x  \right)^{m_j} \left(\tfrac{1}{ 1-x^L} \right)
  + \dots + b_j(L-1) \left( x\partial_x  \right)^{m_j}
  \left( \tfrac{x^{L-1} }{ 1-x^L}  \right) \right).
\end{align*}
We now note that each rational function $\frac{x^j }{ 1-x^L}$ on
the right-hand of the last equation is an eigenfunction of
$U_{L+1}$, which follows easily from its Taylor series. By
Lemma~\ref{IteratedEF}, the same statement holds for $\left(
x\partial_x  \right)^{m_j}  \left( \frac{x^{L-1} }{ 1-x^L}
\right)$. In conclusion,  we have expressed $f$ as a finite linear
combination of real eigenfunctions of the operator $U_{L+1}$.
\end{proof}

The vector space generated by all of the eigenfunctions given in this
theorem is infinite dimensional.  It is natural to ask how we can
decompose it into finite dimensional vector spaces, so that we can do
analysis on each finite dimensional piece with greater ease.
When we fix the weight $\kappa$ and the level $L$ of admissible
eigenfunctions, we obtain a finite-dimensional vector space
(Theorem~\ref{DirectSum} below) of eigenfunctions.  We note that this
grading is quite natural, given the structure theorem.  It also
plays an analogous role to the grading of the finite-dimensional
vector spaces of cusp forms and Eisenstein series that arise in
automorphic forms.

We now define the relevant notions that are used in the
afore-mentioned grading.

\begin{definition}
We denote by
\[ \V_{\kappa, L}(U_p) \]
the vector space of all real rational functions with fixed weight
$\kappa$ and fixed level $L$, that are eigenfunctions of the Hecke
operator $U_p$.  Given a set of integers $S = \{ p_1, \dots, p_n
\}$, we let
\[ \V_{\kappa, L}(S) = \V_{\kappa, L}(U_{p_1}, \dots, U_{p_n}) \]
denote the vector space over $\R$ of all real rational functions
with fixed weight $\kappa$ and fixed level $L$ that are
simultaneous eigenfunctions of the collection of Hecke operators
$U_{p_1}, \dots, U_{p_n}$.

We remark that when $\{ p_1, p_2, p_3, \dots,  p_n \} $ is the set
of all integers between $2$ and $L$ inclusively, the corresponding
vector space of simultaneous eigenfunctions for $U_2, U_3, \dots,
U_L$ is 1-dimensional, generated by the function
\[ \sum_{n=0}^\infty  \chi(n) n^{\kappa -1} x^n,  \]
as we saw in Theorem~\ref{Character} of the previous section (with
$\chi$ being the real character mod $L$).  We further define
\begin{align*}
\S_{\kappa, L}(S) &= \S_{\kappa, L} (U_{p_1}, \dots, U_{p_n} ) =
\{ f \in  \V_{\kappa, L}(S)  \st a_0 = 0 \}.
\end{align*}
\end{definition}

It is clear that $ \S_{\kappa, L}(S) $ is a vector space over
$\mathbb{R}$, and is $U_{p_j}$-invariant for each $p_j \in S$.

\begin{theorem}\label{DirectSum}
For a fixed weight $\kappa$ and fixed level $L$, $\V_{\kappa,
L}(S)$ is a finite-dimensional vector space. Considered as a
vector space over $\mathbb{C}$, it has the basis
\[ \{ f_{\chi_1}, f_{\chi_2}, \dots,  f_{\chi_{\phi(L)}} \}, \]
where $\phi(L)$ is the Euler $\phi$-function of $L$, and where
\[  f_{\chi}(x) = \sum_{n=0}^\infty \chi(n)n^{\kappa -1} x^n.   \]
\end{theorem}
\begin{proof}
We make the easy observation that for each fixed denominator
$B(x)$ of an eigenfunction there are at most finitely many
possible numerators, each numerator being an eigenfunction of the
corresponding matrix defined in Section $3$. Note that the degree
of $B(x)$ must be less than or equal to $\kappa L$, by the
structure theorem. Since there are at most finitely many possible
denominators of degree $\leq \kappa L$ whose roots are $L$'th
roots of unity, we conclude that there are at most finitely many
linearly independent eigenfunctions of level $L$ and weight
$\kappa$.  The second statement concerning the basis is tantamount
to doing Fourier analysis on the finite group
$\mathbb{Z}/L\mathbb{Z}$ (see Knapp~\cite{MR93j:11032}, for
example), from which we know that we can expand every periodic
function into a complex linear combination of the $\phi(L)$
characters $\chi$ mod $L$.
\end{proof}

The dimensions of the vector spaces defined here offer challenging
combinatorial problems. Indeed, it is not clear how to compute
dim$(\V_{1,L}(U_2) )$ even in the case when $L$ is prime (and
involves the Artin conjecture for primitive roots mod $L$).

\begin{example}
We note that the space $\V_{1,L}$ always has the eigenfunction
\[
  f(x) = \frac{2+x}{1 - 2\cos(\frac{2\pi}{L}) x +x^2},
\]
of weight $1$ and level $L$.  The Taylor coefficients of $f =
\sum_{n=0}^\infty a_n x^n$ are given by
$ a_n = e^\frac{2\pi i n}{L} + e^\frac{-2\pi i n}{L}
= 2\cos(\frac{2\pi n}{L})$.
Equivalently, the Taylor coefficients satisfy the linear recurrence
\[
a_n =  2\cos(\tfrac{2\pi}{L})a_{n-1} - a_{n-2}.
\]
Indeed, we have $U_{L-1}f = f$. We single out this class of
eigenfunctions for being eigenfunctions of $U_{L-1}$, but of no
other Hecke operator $U_p$ with $p < L$. To wit,
\[
 a_{pn}= e^\frac{2\pi i pn}{L} +
 e^\frac{-2\pi i pn}{L}= \lambda_p a_n
\]
for some $\lambda_p$, only when $p = L-1$ (assuming $p < L$).
\end{example}

%%%%%%%%%%%%%%%%%%%%%%%%%%%%%%%%%%%%%%%%%%%%%%%%%%%%%%%%%%%%%%%%%%%%
\section{Simultaneous eigenfunctions}\label{sec_Simultaneous}
%%%%%%%%%%%%%%%%%%%%%%%%%%%%%%%%%%%%%%%%%%%%%%%%%%%%%%%%%%%%%%%%%%%%

Consider the algebra $\mathfrak{H}=\{U_p \st p\in\N, \; p\ge 2\}$
of all Hecke operators acting on rational functions. We are
interested in the intersection of the spectra and in the set of
common eigenfunctions.  Recall that $\lambda=1$ belongs to
$\spec(U_p)$ for every $p$, and $\frac1{1-x}$ is a common
eigenfunction for the full algebra $\mathfrak{H}$.  In this
section we give a precise description of all possible common
eigenvalues and eigenfunctions for the whole algebra
$\mathfrak{H}$.

We first show that the simultaneous spectrum of any two Hecke
operators $U_m$ and $U_n$ is trivialized, if $m$ and $n$ are
relatively prime. However, it turns out that simultaneous
eigenspaces are in general non-trivial.

%Indeed it is the consideration of a function $f$ that is
%simultaneously an eigenfunction of $U_p$ for $p =2, 3, \dots, L$ that
%gives us a character mod $L$.

\begin{lemma}
For any two relatively prime integers $m$ and $n$, we have
\[ \spec(U_m) \cap \spec(U_n)= \{0,\pm 1\}. \]
\end{lemma}
\begin{proof}
We already know that $0$, $1$ and $-1$ are always contained in the
spectrum of $U_p$ for every $p$,  from
Lemma~\ref{IntegerSpectrum}. Thus we only need to show the other
inclusion.  In fact, if $\lambda\in\spec(U_m)\cap\spec(U_n)$, then
either $\lambda = 0$ or Theorem~\ref{EValueStructure} implies
$\lambda=\pm m^k$ and $\lambda=\pm n^\ell$ for some $k,\ell\in\N$.
But this implies $k=\ell=0$ since $(m,n)=1$, that is, $\lambda=\pm 1$.
\end{proof}

%%%%%%%%%%%%%%%%%%%%%%%%%%%%%%%%%%%%%%%%%%%%%%%%%%%%%%%%%%%%%%%%%%%%
\begin{lemma}\label{f=zero}
Let $f$ be such that $U_p f=\lambda f$ for some $p$, $\lambda\not=0$. 
If $f(x)=x^m \tilde f(x)$ for some positive integer $m$ and
some $\tilde f$ with $\tilde f(0)\not=0$, then $p \nmid m$.
\end{lemma}
\begin{proof}
Using
\eqref{f-identity} we get
\begin{align*}
\lambda\; x^{pm} \tilde f(x^p)&= \frac1p \sum_{j=0}^{p-1}
 e^{\frac{2\pi i m}{p}j} x^m\tilde f(e^{\frac{2\pi i}{p}j} x)
\intertext{which implies}
\lambda\; x^{pm-m} \tilde f(x^p)&= \frac1p \sum_{j=0}^{p-1}
 e^{\frac{2\pi i m}{p}j} \tilde f(e^{\frac{2\pi i}{p}j} x).
\end{align*}
Evaluating at $x=0$ gives
\[ 0=\frac1p \sum_{j=0}^{p-1} e^{\frac{2\pi i m}{p}j}\tilde f(0) \]
and therefore $p \nmid m$.
\end{proof}

%%%%%%%%%%%%%%%%%%%%%%%%%%%%%%%%%%%%%%%%%%%%%%%%%%%%%%%%%%%%%%%%%%%%
The following theorem completely describes the set of simultaneous
eigenfunctions for the algebra $\mathfrak{H}$ of all Hecke operators on
rational functions. Among others, our description reveals the importance
of the operator $x\partial_x$.

\begin{theorem}[Simultaneous Eigenfunctions]\label{simultaneous}
Let $f$ be a simultaneous eigenfunction of $\mathfrak{H}$ such
that $f$ is not in the kernel of $U_p$ for any $p$. That is, let
$f$ be a rational function with the property that for
\textit{every} $p$ there is a $\lambda_p\not=0$ such that $U_p
f=\lambda_p f$. Then $L=1$, and
\[ f(x) = C(x\partial_x)^k \left(\frac1{1-x}\right)\]
for some $k\in\N$ and $C\in\C$. Consequently, $\lambda_p=p^k$.
\end{theorem}
\begin{proof}
We consider the two cases $f(0)\not=0$ and $f(0)=0$ separately. If
$f(0)\not=0$, then we plug in $x=0$ into \eqref{f-identity} and
get $\lambda_p=1$, so $U_p f= f$ for every $p$. If we write
$f(x)=\sum_{n=0}^\infty a_n x^n$, then we must have $a_{p}=a_1$
for every $p$, so
\begin{align*}
 f(x) &= a_0 + a_1 x\sum_{n=0}^\infty x^n \\
 &= a_0 +\frac{a_1 x}{1-x} = \frac{a_0+(a_1-a_0)x}{1-x}
\end{align*}
which implies $a_1=a_0$, in other words, $f(x)=a_0(\frac1{1-x})$ as claimed.

Suppose now $f(0)=0$ and write $f(x)=x^m \tilde f(x)$ with $\tilde
f(0)\not=0$. By Lemma~\ref{f=zero}, $f$ cannot be an eigenfunction
of $U_{m}$. Since $f$ is assumed to be a common eigenfunction of
$\mathfrak{H}$, it follows that $m=1$. The structure theorem shows
that if $U_p f=\lambda_p f$ with $\lambda_p\not=0$, then
$\lambda_p = \chi_f(p)\, p^k$ for some integer $k$ and
$\chi_f(p)=\lambda_p/|\lambda_p|$. Write $f(x)=\sum_{n=1}^\infty
a_n x^n$ with $a_1\not=0$. Then $a_{p}=\chi_f(p)\, p^k a_1$ holds
for every $p$, thus
\begin{align*}
 f(x) &= \sum_{n=1}^\infty \chi_f(n)\, n^k a_1 x^n \\
 &= a_1 \sum_{n=1}^\infty n^k x^n \\
 &= a_1 (x\partial_x)^k \left(\frac1{1-x}\right)
\end{align*}
since $\chi_f(n) = 1$ for every $n$ by Corollary~\ref{Cor4-3}(iii).
\end{proof}

%%%%%%%%%%%%%%%%%%%%%%%%%%%%%%%%%%%%%%%%%%%%%%%%%%%%%%%%%%%%%%%%%%%%
We denote by $\V$ the vector space over $\C$ spanned by the functions
\begin{equation}\label{SEigenfunctions}
 \phi_k(x) = (x\partial_x)^k \left(\frac1{1-x}\right)
\end{equation}
for $k\in\N$. That is,
\begin{equation}\label{SEigenspace}
 \V = \LinSpan_{\C}\{\phi_0,\phi_1,\phi_2,\dots\}.
\end{equation}

\begin{remark}
Notice that although every $\phi_k$ is a simultaneous
eigenfunction of $\mathfrak{H}$, the sum of two of them
$\phi_i+\phi_j$ is not, simply because their weights are different
($i\not=j$). However, if $f=\sum c_k \phi_k\in\V$, then $U_p f =
\sum c_k U_p \phi_k = \sum c_k p^k \phi_k \in \V$. Thus the space
$\V$ is $U_p$-invariant for every $p$.
\end{remark}

%%%%%%%%%%%%%%%%%%%%%%%%%%%%%%%%%%%%%%%%%%%%%%%%%%%%%%%%%%%%%%%%%%%%
Every function $\phi_k$ can obviously be written as
\[ \phi_k(x) = \frac{A_k(x)}{(1-x)^{k+1}} \]
where
\[ A_k(x) = (1-x)^{k+1} \sum_{n=1}^\infty n^k x^n. \]
It is easy to check the identity
\begin{equation}\label{EulerianPoly}
 A(x) = \sum_{\ell=0}^k S(k,\ell)\,\ell!\, x^\ell (1-x)^{k-\ell}
\end{equation}
where $S(k,\ell)$ are the well-known Stirling numbers of the second
kind.

The polynomials $A_k$ are known as Eulerian polynomials,
cf.~\cite{MR57:124}.
Here are the first few:
\begin{align*}
A_1(x) &= x, \quad A_2(x) = x+x^2, \\
A_3(x) &= x+4x^2+x^3, \\
A_4(x) &= x+11x^2+11x^3+x^4, \\
A_5(x) &= x+26x^2+66x^3+26x^4+x^5, \\
A_6(x) &= x+57x^2+302x^3+302x^4+57x^5 +x^6. \\
\end{align*}

%%%%%%%%%%%%%%%%%%%%%%%%%%%%%%%%%%%%%%%%%%%%%%%%%%%%%%%%%%%%%%%%%%%%
\begin{lemma}
The family of functions $\{\phi_k \st k\in\N\}$ from
\eqref{SEigenfunctions} is a linearly independent system. In
particular, $\V$ is an infinite dimensional vector space.
\end{lemma}
\begin{proof}
For any given integer $n$ we will prove that the functions $\phi_0,
\dots,\phi_n$ are linearly independent. For arbitrary constants
$c_0,\dots,c_n$, we have
\begin{align*}
 \sum_{k=0}^n c_k\phi_k(x)
 &= \sum_{k=0}^n c_k \frac{A_k(x)}{(1-x)^{k+1}} \\
 &= \sum_{k=0}^n \sum_{\ell=0}^k c_k \alpha_\ell
\frac{x^\ell}{(1-x)^{\ell+1}}
\intertext{(where $\alpha_\ell=S(k,\ell)\,\ell!$ from
 \eqref{EulerianPoly}) }
 &= \sum_{\ell=0}^n \sum_{k=\ell}^n c_k \alpha_\ell
\frac{x^\ell}{(1-x)^{\ell+1}} \\
 &= \frac{1}{(1-x)^{n+1}} \sum_{\ell=0}^n \alpha_\ell
\bigg(\sum_{k=\ell}^n c_k\bigg) x^\ell (1-x)^{n-\ell}.
\end{align*}
If $\sum_{k=0}^n c_k\phi_k(x)=0$, then
$\alpha_\ell \sum_{k=\ell}^n c_k = 0$ for every $\ell$; that is,
\begin{equation*}
\begin{pmatrix}
 \alpha_0 & \alpha_0 & \alpha_0 & \dots & \alpha_0 \\
 0 & \alpha_1 & \alpha_1 & \dots & \alpha_1 \\
 \vdots & \ddots & \alpha_2 & \dots & \alpha_2 \\
 \vdots & & \ddots & \ddots & \vdots \\
 0 & \dots & \dots & 0 & \alpha_n
\end{pmatrix}
\begin{pmatrix}
 c_0 \\ c_1 \\ c_2 \\ \vdots \\ c_n
\end{pmatrix}
= \mathbf{0}.
\end{equation*}
Thus $c_0=c_1=\cdots = c_n=0$ since the $\alpha_\ell$ are all different
from zero.
\end{proof}
\begin{lemma}\label{VEigenfunction}
If $f\in\V$ is an eigenfunction for some $U_p$, then it is a
simultaneous eigenfunction of the whole Hecke algebra
$\mathfrak{H}$.
\end{lemma}
\begin{proof}
Let $0\not= f =\sum c_j \phi_j\in\V$ satisfy $U_pf= \chi_f(p) p^k
f$ for some $p$ and some $k$. We also have
\begin{align*}
  U_pf = \sum c_j U_p \phi_j = \sum c_j p^j \phi_j
\intertext{so that we get the identity}
 \sum c_j (p^j - \chi_f(p) p^k) \phi_j = 0.
\end{align*}
The linear independence of the $\phi_j$ implies $c_j=0$ for every 
$j\not=k$ and therefore $f= c_k \phi_k$, a simultaneous eigenfunction 
of $\mathfrak{H}$.
\end{proof}

%%%%%%%%%%%%%%%%%%%%%%%%%%%%%%%%%%%%%%%%%%%%%%%%%%%%%%%%%%%%%%%%%%%%
\section{A first application: tensor products of Hecke operators \\
and the Riemann zeta function} \label{sec_TensorProduct}
%%%%%%%%%%%%%%%%%%%%%%%%%%%%%%%%%%%%%%%%%%%%%%%%%%%%%%%%%%%%%%%%%%%%

In this section we consider tensor products of Hecke operators, and we
study their spectral properties. Let $\V$ be the vector space
introduced in \eqref{SEigenspace} and let
\[ \H^n = \underbrace{\V\otimes \cdots\otimes \V}_{n \text{ times}}. \]
For any finite set $S=\{p_1,\dots,p_n\}$ of prime numbers with $p_j\not=
p_k$ if $j\not=k$, we define the operator
\[ \U_S=U_{p_1}\otimes\cdots\otimes U_{p_n}:
   \H^n \to \rational(x_1,\dots,x_n) \]
by
\begin{equation}\label{TensorOperator}
 \U_S(f_1\otimes\cdots\otimes f_n)(x_1,\dots,x_n) =
    (U_{p_1} f_1)(x_1) \cdots (U_{p_1} f_n)(x_n),
\end{equation}
where every $U_{p_k}$ is a Hecke operator and the multiplication
on the right-hand side of \eqref{TensorOperator} is the usual
multiplication of rational functions.

%%%%%%%%%%%%%%%%%%%%%%%%%%%%%%%%%%%%%%%%%%%%%%%%%%%%%%%%%%%%%%%%%%%%
\begin{lemma}\label{TensorEigenfunction}
An element $f_1\otimes\cdots\otimes f_n\in\H^n$ is an
eigenfunction of $\U_S$ if and only if each $f_j$ is a
simultaneous eigenfunction of the Hecke algebra $\mathfrak{H}$.
\end{lemma}
\begin{proof}
First, let each $f_j$ be a simultaneous eigenfunction of $\mathfrak{H}$.
Thus $\chi_{f_j}=1$ (the identity character) and we have
\begin{equation}\label{TS-Spectrum}
\begin{split}
 \U_S(f_1\otimes\cdots\otimes f_n)
 &= (U_{p_1} f_1)\dots (U_{p_1} f_n) \\
 &= (p_1^{k_1}f_1) \dots (p_n^{k_n} f_n) \\
 &= (p_1^{k_1}\dots p_n^{k_n}) f_1\otimes\cdots\otimes f_n.
\end{split}
\end{equation}
Therefore $f_1\otimes\cdots\otimes f_n$ is an eigenfunction of
$\U_S$ with eigenvalue $p_1^{k_1}\dots p_n^{k_n}$ for some
integers $k_1,\dots,k_n$.

To prove the converse we now suppose that
\[ \U_S(f_1\otimes\cdots\otimes f_n)
   = \lambda f_1\otimes\cdots\otimes f_n \]
for some $\lambda$, and some $f_1,\dots,f_n\in\V$ such that no $f_j$ is
the zero function. Then we have
\[ \lambda = \left(\frac{U_{p_1}f_1(x_1)}{f_1(x_1)}\right) \cdots
   \left(\frac{U_{p_n}f_n(x_n)}{f_n(x_n)}\right) \]
which implies that every factor
$\left(\frac{U_{p_j}f_j(x_j)}{f_j(x_j)}\right)$ must be constant.
Thus, for every $j$, $f_j$ is an eigenfunction of $U_{p_j}$ and
consequently a simultaneous eigenfunction of $\mathfrak{H}$ by
Lemma~\ref{VEigenfunction}.
\end{proof}

%%%%%%%%%%%%%%%%%%%%%%%%%%%%%%%%%%%%%%%%%%%%%%%%%%%%%%%%%%%%%%%%%%%%
\begin{corollary}\label{Multiplicity}
The spectrum of $\U_S$ is the set
\[ \{p_1^{k_1}\cdots p_n^{k_n} \st k_j\in\N
   \text{ for every } j\}, \]
where each eigenvalue occurs with multiplicity $1$.
\end{corollary}
\begin{proof}
 From \eqref{TS-Spectrum} it is clear that every integer
$p_1^{k_1}\cdots p_n^{k_n}$ lies in the spectrum of $\U_S$. On the
other hand, it also follows from the proof of the previous theorem
that, if $\lambda$ is an eigenvalue of $\U_S$, then
$\lambda=\lambda_1\cdots\lambda_n$ with $U_{p_j}f_j=\lambda_j f_j$
for some $f_j\in\V$. Lemma~\ref{VEigenfunction} then implies that
$f_j$ is a simultaneous eigenfunction, so $f_j = c_j \phi_{k_j}$
for some integer $k_j$ and some constant $c_j$. Thus  $\lambda =
\lambda_1 \cdots \lambda_n = p_1^{k_1} \cdots  p_n^{k_n}$.

To prove the multiplicity $1$ statement, let $g_1\otimes\cdots\otimes
g_n$ be another eigenfunction for $\lambda$. Then we similarly get $g_j
= c_j' \phi_{\ell_j}$ and so $\lambda=p_1^{\ell_1}\cdots p_n^{\ell_n}$.
Finally, by the unique factorization theorem for $\mathbb{Z}$, it
follows that $k_j=\ell_j$ for every $j$. Hence each $g_j$ is a multiple
of $f_j$ and the multiplicity of $\lambda$ is $1$.
\end{proof}

%%%%%%%%%%%%%%%%%%%%%%%%%%%%%%%%%%%%%%%%%%%%%%%%%%%%%%%%%%%%%%%%%%%%
\begin{definition}
For an operator $A$ with nonnegative discrete spectrum, let
\begin{equation*}
 \zeta_{A}(s) = \sum_{\lambda\in\spec_+(A)}
 \frac{m(\lambda)}{\lambda^s},
\end{equation*}
where $\spec_+(A)=\spec(A)\setminus \{0\}$ is the set of positive
eigenvalues of $A$ and $m(\lambda)$ is their multiplicity.
\end{definition}

%%%%%%%%%%%%%%%%%%%%%%%%%%%%%%%%%%%%%%%%%%%%%%%%%%%%%%%%%%%%%%%%%%%%
\begin{example}
For a Hecke operator $U_p$ acting on $\V$ we have
\begin{equation*}
 \zeta_{U_p}(s) = \sum_{j=0}^\infty \frac1{p^{js}}
 = \frac1{1-p^{-s}}.
\end{equation*}
\end{example}

%%%%%%%%%%%%%%%%%%%%%%%%%%%%%%%%%%%%%%%%%%%%%%%%%%%%%%%%%%%%%%%%%%%%
The following theorem is an interesting application of the
spectral properties of the Hecke operators acting on the vector
space $\V$ spanned by the simultaneous eigenfunctions of
$\mathfrak{H}$.  More precisely, the spectrum of $\U_S$ forges a
natural link to the Riemann zeta function.

\begin{theorem}[Euler product]
Let $S=\{p_1,\dots,p_n\}$ be a set of prime numbers. Then
\begin{equation*}
 \zeta_{\U_S}(s) = \zeta_{U_{p_1}}(s)\cdots \zeta_{U_{p_n}}(s)
 =\prod_{p\in S} \frac1{1-p^{-s}}.
\end{equation*}
In other words, the zeta function of the operator $\U_S$ acting on
$\H^n$ $($cf. \eqref{TensorOperator}$)$ is a finite Euler product
of the Riemann zeta function $\zeta(s)$.
\end{theorem}
\begin{proof}
 From the definition of the zeta function of $\U_S$ and because the
multiplicity $m(\lambda)=1$ for all $\lambda\in\spec(\U_S)$ (by
Corollary~\ref{Multiplicity}), we have
\begin{align*}
 \zeta_{\U_S}(s) &= \sum_{\lambda \in \spec(\U_S)} \frac{1}{\lambda^s}\\
 &= \sum_{k_1, \dots, k_n \in \mathbb{N} }
  \frac{1}{(p_1^{k_1} \cdots p_n^{k_n})^s}\\
 &= \sum_{k_1 = 0}^\infty \left(\frac{1}{p_1^s} \right)^{k_1} \cdots
    \sum_{k_n = 0}^\infty  \left(\frac{1}{p_n^s} \right)^{k_n} \\
 &= \prod_{p\in S} \frac1{1-p^{-s}}.
\end{align*}
\end{proof}

Let $\H^\infty$ be the space of products
$\mathbf{f}= f_1\otimes f_2\otimes \cdots$, where $\{f_n\}_{n\in\N}$
is an infinite sequence of rational functions with the following properties:
\begin{enumerate}
\item
There is a finite set $I\subset \N$ such that $f_j\in \V$ for every $j\in I$.
\item
$f_j= 1$ for every $j\in \N\setminus I$.
\end{enumerate}

For $\mathbf{f}\in \H^\infty$ we define the operator $\U$ by
\begin{equation}\label{BigTensor}
  \U\mathbf{f} =
  (U_{q_1}f_{i_1}) \otimes\cdots\otimes (U_{q_{m}}f_{i_m}),
\end{equation}
where $I=\{i_1,\dots,i_m\}$ is the finite set of positive integers
associated to $\mathbf{f}$, and for each $k$ the number
$q_k=p_{i_k}$ is the $i_k$-th prime number. Notice that similarly
to $\U_S$, the operator $\U$ maps tensor products of rational
functions into rational functions in several variables.

Given $\mathbf{f}\in \H^\infty$ with corresponding index set
$\{i_1,\dots,i_m\}$, we let $S_{\mathbf f}$ be the set of prime numbers
$\{p_{i_1},\dots,p_{i_m}\}$ from \eqref{BigTensor}. Clearly,
\[ \U\mathbf{f} =
   \U_{S_{\mathbf f}}(f_{i_1}\otimes\cdots\otimes f_{i_m}), \]
where $\U_{S_{\mathbf f}}$ is the operator from
\eqref{TensorOperator}. Therefore Lemma~\ref{TensorEigenfunction}
and Corollary~\ref{Multiplicity} apply verbatim to the operator
$\U$.

Since $\Z$ is a unique factorization domain,
Corollary~\ref{Multiplicity} implies that the spectrum of $\U$ is
exactly the set of all positive integers and each eigenvalue has
multiplicity one. This leads to the following result.

\begin{theorem}[Riemann zeta function]
The spectral zeta function of the operator $\U$ on $\H^\infty$
satisfies
\begin{equation*}
 \zeta_{\U}(s) = \zeta(s),
\end{equation*}
where $\zeta(s)$ is the Riemann zeta function.
\end{theorem}

%%%%%%%%%%%%%%%%%%%%%%%%%%%%%%%%%%%%%%%%%%%%%%%%%%%%%%%%%%%%%%%%%%%%
\section{A second application: completely multiplicative coefficients}
\label{sec_Application2}
%%%%%%%%%%%%%%%%%%%%%%%%%%%%%%%%%%%%%%%%%%%%%%%%%%%%%%%%%%%%%%%%%%%%
Mordell and Hecke were motivated to study simultaneous eigenforms
in the context of modular forms in order to classify those forms
that have multiplicative coefficients, after the discovery that
the coefficients of the discriminant function $\Delta(\tau)$
indeed have those properties. In the same spirit, we now give a
complete classification of those rational functions that have
completely multiplicative coefficients, since our Hecke operators
by definition have a completely multiplicative action on rational
functions: $U_m U_n f = U_{mn} f$, for all $m, n$.

We use the vector spaces $\V$ and $\V_{\kappa, L}(U_2, \dots,
U_L)$ of the simultaneous eigenfunctions to characterize all
rational functions $f = \sum_{n=0}^\infty a_n x^n \in \rational$
whose coefficients are completely multiplicative - i.e. such that
$a_{mn} = a_m a_n $ for all $m,n$.

\begin{theorem}
A rational function $f = \sum_{n=0}^\infty  a_n x^n \in \rational$
has completely multiplicative coefficients if and only if it has
the following form:
\begin{itemize}
\item[(i)]    If $a_0 =0$, then there exist positive integers
$\kappa$ and $L \geq 1$ such that
\[ f(x) = a_1 \ (x \partial_x)^{\kappa-1}
  \bigg(\frac{\sum_{j=1}^{L-1} \chi_f(j) x^j}{1-x^L} \bigg),
\]\label{specialform}
where $\chi_f$ is the real quadratic character mod $L$.
\item[(ii)]   If $a_0 \not= 0$, then
\[ f(x) = \frac{a_0}{1-x}. \]
\end{itemize}
\end{theorem}
\begin{proof}
In both cases, we observe that if we \emph{fix} $m$, then the
assumption that $a_{nm} = a_n a_m $ holds for all $n$ can be
regarded as telling us that the coefficients $a_k$ give us an
eigenfunction of $U_m$, with eigenvalue $\lambda = a_n$.
Furthermore, this condition is satisfied by \emph{all} $m$, hence
making the function $f = \sum_{n=0}^\infty  a_n x^n$ a
simultaneous eigenfunction of all the Hecke operators.

We first treat the case $L>1$.  Here the
hypotheses of Theorem~\ref{Character}(i) are satisfied by the
remarks made in the previous paragraph concerning simultaneous
eigenfunctions, giving us the eigenfunctions in part (i) of the
theorem.

For the case $L=1$, it is trivially true that
$\chi_f(p)=1$ for all $p$.  The conclusion of Corollary~\ref{Cor4-3}(iii)
now implies that $f$ is not in the kernel of any of the Hecke operators.
Therefore the hypotheses of
Theorem~\ref{simultaneous} are satisfied, giving us the eigenfunctions
in part (i) of the theorem, where in this case $L=1$ and $\chi_f(j)$ is 
the trivial character.
\end{proof}

%%%%%%%%%%%%%%%%%%%%%%%%%%%%%%%%%%%%%%%%%%%%%%%%%%%%%%%%%%%%%%%%%%%%
\section{Appendix: Explicit examples of eigenfunctions} \label{sec_Examples}
%%%%%%%%%%%%%%%%%%%%%%%%%%%%%%%%%%%%%%%%%%%%%%%%%%%%%%%%%%%%%%%%%%%%

For illustration purposes, we now focus on the finite dimensional
vector spaces $\S_{\kappa,L}(U_2)$ of eigenfunctions with constant
term equal to zero ($a_0=0$), in the range $\kappa L\le 6$. This
class is large enough to already exhibit some of the non-trivial
features of the eigenspaces. As noted in
Remark~\ref{DegreeDenominator} the only possible degree for the
denominator of such an eigenfunction is $d=1,\dots,6$.

We recall that $\lambda=1$ is the only possible eigenvalue
associated to an eigenfunction $f$ such that $f(0)\not=0$.
All other
eigenvalues must come from the vector space $\S_{\kappa,L}(U_2)$.
Fortuitously, the case $\lambda=1$ has recently been investigated
in \cite{BLMMS02} for an operator related to $U_2$. We therefore
concentrate in this Appendix on the eigenvalues $\lambda\not=1$
and consequently on $\S_{\kappa,L}(U_2)$.

Since we are  only considering eigenfunctions here that do not lie
in the kernel of $U_2$, let $\lambda\not=0$. Recall that by
Lemma~\ref{f-splitting}, $U_2 f = \lambda f$ if and only if
\begin{equation}\label{f2-identity}
 \lambda f(x^2) = \frac12(f(x) + f(-x)).
\end{equation}

Let $f(x)=A(x)/B(x)$ be an eigenfunction of $U_2$ with
\[ B(x)=1+\alpha_1 x +\cdots + \alpha_d x^d,\; \alpha_d\not=0. \]
In the case $p=2$, the identity \eqref{B-identity} becomes
\begin{equation}\label{B2-identity}
 B(x^2)=B(x)B(-x).
\end{equation}
In other words,
\begin{align*}
1+\alpha_1 x^2 &+ \cdots + \alpha_d x^{2d}\\
 &= (1+\alpha_1 x +\cdots + \alpha_d x^d)
    (1-\alpha_1 x +\cdots + (-1)^d\alpha_d x^d)\\
 &= 1+(2\alpha_2-\alpha_1^2) x^2 +\cdots + (-1)^d\alpha_d^2 x^{2d}
\end{align*}
leading to the formula
\begin{equation}\label{B2-coeff}
 \alpha_\ell =
 \sum_{\substack{j+k=2\ell\\j,k\ge 0}}(-1)^j \alpha_j\alpha_k
 \;\text{ for } \ell=1,\dots,d,
\end{equation}
where $\alpha_0=1$. Moreover, due to the involution property $x^d
B(\frac1x)=(-1)^d B(x)$ from Theorem~\ref{B-structure}, we have
$\alpha_j=(-1)^d\alpha_{d-j}$ for every $j$. This relation
together with \eqref{B2-coeff} allow us to find the possible
coefficients for the denominator of an eigenfunction of $U_2$. In
order to find all possible numerators with a given denominator
$B(x)$, we further use the identity
\begin{equation}\label{A2-identity}
2\lambda A(x^2)= A(x)B(-x) +A(-x)B(x)
\end{equation}
which follows from \eqref{f2-identity} and \eqref{B2-identity}.

%%%%%%%%%%%%%%%%%%%%%%%%%%%%%%%%%%%%%%%%%%%%%%%%%%%%%%%%%%%%%%%%%%%%
By writing $A(x)=\displaystyle\sum_{j=1}^{d-1} c_j x^j$ (note that 
$f(0) = 0$ implies $A(0)=0$), the identity \eqref{A2-identity} becomes
\begin{align*}
 \lambda \sum_{j=1}^{d-1} c_j x^{2j}
 &= \sum_{j=1}^{d-1}\left(\sum_{k=1}^{2j}
 (-1)^{k} \alpha_{2j-k} c_k\right) x^{2j}
\intertext{implying}
 \lambda c_j &= \sum_{k=1}^{2j} (-1)^{k} \alpha_{2j-k} c_k
\end{align*}
for every $j=1,\dots,d-1$.
Setting $c_j:=0$ for $j>d-1$ and $\alpha_k:=0$ for $k>d$ we get
that $\lambda$ is an eigenvalue of the matrix
$\mathfrak{B} = ((-1)^{k} \alpha_{2j-k})_{j,k=1,\dots,d-1}$.
We now give the explicit form of $\mathfrak{B}$ for various values of $d$.
\begin{itemize}
\item $d=3$:
\begin{equation*}
\mathfrak{B}=
\begin{pmatrix}
-\alpha_1 & 1  \\
1 & -\alpha_1
\end{pmatrix}
\end{equation*}
since $\alpha_3=-1$ and $\alpha_1+\alpha_2=0$.
\item $d=4$:
\begin{equation*}
\mathfrak{B}=
\begin{pmatrix}
-\alpha_1 & 1 & 0 \\
-\alpha_1 & \alpha_2 & -\alpha_1 \\
0 & 1 & -\alpha_1
\end{pmatrix}
\end{equation*}
since $\alpha_4=1$ and $\alpha_1-\alpha_3=0$.
\item $d=5$:
\begin{equation*}
\mathfrak{B}=
\begin{pmatrix}
-\alpha_1 & 1 & 0 & 0\\
\alpha_2 & \alpha_2 & -\alpha_1 & 1\\
1 & -\alpha_1 & \alpha_2 & \alpha_2\\
0 & 0 & 1 & -\alpha_1
\end{pmatrix}
\end{equation*}
since $\alpha_5=-1$ and $\alpha_1+\alpha_4=\alpha_2+\alpha_3=0$.
\item $d=6$:
\begin{equation*}
\mathfrak{B}=
\begin{pmatrix}
-\alpha_1 & 1 & 0 & 0 & 0\\
-\alpha_3 & \alpha_2 & -\alpha_1 & 1 & 0\\
-\alpha_1 & \alpha_2 & -\alpha_3 & \alpha_2 & -\alpha_1\\
0 & 1 & -\alpha_1 & \alpha_2 & -\alpha_3\\
0 & 0 & 0 & 1 & -\alpha_1
\end{pmatrix}
\end{equation*}
since $\alpha_6=1$ and $\alpha_1-\alpha_5=\alpha_2-\alpha_4=0$.
\end{itemize}

%%%%%%%%%%%%%%%%%%%%%%%%%%%%%%%%%%%%%%%%%%%%%%%%%%%%%%%%%%%%%%%%%%%%
We now list a basis of eigenfunctions of $U_2$ together with their
corresponding eigenvalues $\lambda$ for every given degree $d$.
Using the fact that $|\lambda|=2^{\kappa-1}$ (for $U_2$) reduces the
computations for the admissible denominators $B(x)$.    
\begin{itemize}
\item $d=2$:
\begin{alignat*}{2}
f_{2,1}(x) &=\frac{x}{(1-x)^2}, & \qquad \lambda &= 2, \\
f_{2,2}(x) &=\frac{x}{1+x+x^2}, & \lambda &= -1.
\end{alignat*}
\item $d=3$:
\begin{alignat*}{2}
f_{3,1}(x) &=\frac{x+x^2}{(1-x)^3}, & \qquad \lambda &= 4, \\
f_{3,2}(x) &=\frac{x+x^2}{1-x^3}, & \lambda &= 1.
\end{alignat*}
\item $d=4$:
\begin{alignat*}{2}
f_{4,1}(x) &=\frac{x+4x^2+x^3}{(1-x)^4}, & \qquad \lambda &= 8, \\
f_{4,2}(x) &=\frac{x-x^3}{(1+x+x^2)^2}, & \lambda &= -2, \\
f_{4,3}(x) &=\frac{x+4x^2+x^3}{(1+x+x^2)^2}, & \lambda &= 2, \\
f_{4,4}(x) &=\frac{x-x^3}{1+x+x^2+x^3+x^4}, & \lambda &= -1.
\end{alignat*}
\item $d=5$:
\begin{alignat*}{2}
f_{5,1}(x) &=\frac{x+11x^2+11x^3+x^4}{(1-x)^5}, & \qquad \lambda &= 16, \\
f_{5,2}(x) &=\frac{x+x^2+x^3+x^4}{1-x^5}, & \lambda &= 1.
\end{alignat*}
\item $d=6$:
\begin{alignat*}{2}
f_{6,1}(x) &=\frac{x+26x^2+66x^3+26x^4+x^5}{(1-x)^6}, & \qquad \lambda &= 32, \\
f_{6,2}(x) &=\frac{x-x^2-6x^3-x^4+x^5}{(1+x+x^2)^3}, & \lambda &= -4, \\
f_{6,3}(x) &=\frac{-x-7x^2+7x^4+x^5}{(1+x+x^2)^3}, & \lambda &= 4, \\
f_{6,4}(x) &=\frac{x^3}{(1-x)^2(1+x+x^2)^2}, & \lambda &= 2, \\
f_{6,5}(x) &=\frac{x+2x^2+2x^4+x^5}{(1-x)^2(1+x+x^2)^2}, & \lambda &= 2, \\
f_{6,6}(x) &=\frac{x^3}{1+x^3+x^6}, & \lambda &= -1, \\
f_{6,7}(x) &=\frac{2x+2x^2+x^3}{(1+x+x^2)(1+x+x^2+x^3+x^4)}, & \lambda &= -1, \\
f_{6,8}(x) &=\frac{x+3x^2-3x^4-x^5}{(1+x+x^2)(1+x+x^2+x^3+x^4)}, 
           & \lambda &= 1, \\
f_{6,9}(x) &=\frac{x+2x^2+x^3+2x^4+x^5}{1+x+x^2+x^3+x^4+x^5+x^6}, 
           & \lambda &= 1.
\end{alignat*}
\end{itemize}

We now reorganize the eigenfunctions according to their level $L$,
and give an explicit basis for each $\S_{\kappa, L}(U_2)$ in the
range $\kappa L \leq 6$. Notice that by Corollary~\ref{Cor4-2} we
must have $(p,L)=(2,L)=1$. Thus all of the vector spaces
$\S_{\kappa,2m}(U_2)$ are empty (for all $m$).
\begin{itemize}
\item $L=1$:
\begin{align*}
 \S_{\kappa,1}(U_2) = \LinSpan\{f_{\kappa,1}\}
\end{align*}
for $\kappa=1,\dots,6$.
\item $L=3$:
\begin{align*}
 \S_{1,3}(U_2) = \LinSpan\{f_{2,2},f_{3,2}\}, \quad
 \S_{2,3}(U_2) = \LinSpan\{f_{4,2},f_{4,3},f_{6,4}\}
\end{align*}
\item $L=5$:
\begin{equation*}
 \S_{1,5}(U_2) = \LinSpan\{f_{4,4}, f_{5,2}\}
\end{equation*}
\end{itemize}
Note that the eigenfunctions that do not explicitly appear 
above belong to the spaces $f_{6,2}, f_{6,3} \in \S_{3,3}(U_2)$, 
$f_{6,5}\in \S_{2,3}(U_2)$, $f_{6,6}\in \S_{1,9}(U_2)$,
$f_{6,7}, f_{6,8}\in \S_{1,15}(U_2)$, and
$f_{6,9}\in \S_{1,7}(U_2)$.
We observe that there are some interesting relations among the
eigenfunctions, for example $f_{6,5}=f_{4,3}-6f_{6,4}$.

In addition, we give a basis for $\S_{1,7}(U_2)$
by using the eigenfunctions from Example~\ref{Level7}:
\begin{equation*}
\S_{1,7}(U_2) = \LinSpan\Big\{ \tfrac{x+x^2+x^4}{1-x^7},
\tfrac{x^3+x^5+x^6}{1-x^7}\Big\}.
\end{equation*}

Throughout this section we have made extensive use of MAPLE~6
to compute the eigenfunctions of $U_2$ by finding the
eigenvectors of the corresponding matrix $\mathfrak{B}$ associated
to the given degree.

%%%%%%%%%%%%%%%%%%%%%%%%%%%%%%%%%%%%%%%%%%%%%%%%%%%%%%%%%%%%%%%%%%%%
\bibliographystyle{amsplain}
\providecommand{\bysame}{\leavevmode\hbox to3em{\hrulefill}\thinspace}
\providecommand{\MR}{\relax\ifhmode\unskip\space\fi MR }

\end{document}